\documentclass[a4paper,english]{article}

\usepackage{a4wide}

\usepackage{amsmath,amssymb,amsthm}

\usepackage[english]{babel}

\usepackage[colorlinks=true,linkcolor=blue]{hyperref}

\usepackage{mathrsfs}

\usepackage{siunitx}

\usepackage{booktabs}

\usepackage{stmaryrd}

\usepackage{enumerate}
\usepackage{bm}

\numberwithin{equation}{section}

\usepackage{tikz}
\usetikzlibrary{calc}
\usetikzlibrary{arrows.meta}
\usetikzlibrary{positioning}

\definecolor{myblue}{HTML}{1e77b4}
\definecolor{myorange}{HTML}{ff7f0f}

\definecolor{mycolor}{rgb}{0.122, 0.435, 0.698}

\theoremstyle{definition}

\theoremstyle{definition}

\usepackage{physics}
\usepackage{stmaryrd}

\usepackage{subfig}

\usepackage{pgfplots}
\pgfplotsset{width=.45\textwidth,compat=newest}
\newcommand{\minmod}{\operatorname{minmod}}
\newcommand{\dxi}{\operatorname{d}\xi}

\newcommand{\diag}{\operatorname{diag}}

\newcommand{\sign}{\operatorname{sign}}

\newcommand{\bu}{{u}}

\newcommand{\bw}{{w}}

\title{Numerical investigations into a model of partially incompressible two-phase flow in pipes}

\author{\textsc{Nils Henrik Risebro}\thanks{Department of Mathematics, University of Oslo, Norway (nilshr@math.uio.no, adrianru@math.uio.no)\newline
This project has received funding from the European Union’s Framework Programme for Research and Innovation Horizon 2020 (2014-2020) under the Marie Sk{\l}odowska-Curie Grant Agreement No. 642768.} \and \textsc{Adrian M. Ruf}\footnotemark[1]}
\date{}

\begin{document}

\maketitle

\begin{abstract}
  We consider a model for flow of liquid and gas in a pipe. We assume that the gas is \emph{ideal} and that the liquid is \emph{incompressible}. Under this assumption the resulting system of equations, expressing conservation of mass and momentum, splits into two subsystems such that the gas flow is independent of the liquid flow, and the liquid flow is described by a conservation law parameterized by the mass fraction of gas. When solving these equations numerically, we propose to stagger the gas and liquid variables with respect to each other. The advantage of this is that in finite volume methods one can use numerical flux functions designed for $2\times 2$ systems of hyperbolic conservation laws to solve both the gas flow and the liquid flow, rather than a much more complicated numerical flux for the whole $4\times 4$ system. We test this using the Roe numerical flux for both subsystems, and compare the results with results produced by using the second-order Nessyahu--Tadmor scheme for the second subsystem.
\end{abstract}

% \tableofcontents

%\setlength{\parindent}{0em}

\section{Introduction}\label{sec:intro}
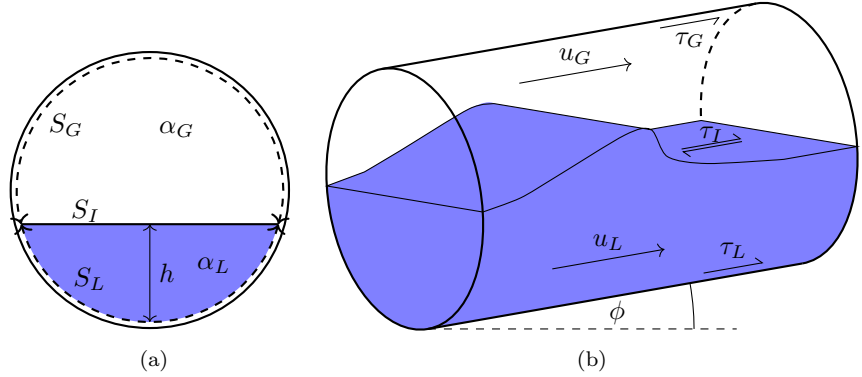
\begin{figure}
  \begin{center}
    \subfloat[]{
      \begin{tikzpicture}
        \fill[blue!50!white] (-1.69,-0.4529) -- (1.69,-0.4529) arc
        (-15:-165:1.75); \draw[thick] (-1.69,-0.4529) --
        (1.69,-0.4529) node[near start,above=-1mm] {$S_I$};
        \draw[thick] (0,0) circle (1.83); \draw[<->,dashed,thick]
        (1.69,-0.4529) arc (-15:195:1.75) node[near end,below
        right=-1mm] {$S_G$}; \draw[<->,dashed,thick] (1.69,-0.4529)
        arc (-15:-165:1.75) node[near end,above right=-1mm] {$S_L$};
        \draw[<->] (0,-1.75) -- (0,-0.4529) node[midway,right] {$h$};
        \node[right] at (0,0.8) {$\alpha_G$}; \node[right] at (0.5,-1)
        {$\alpha_L$};
      \end{tikzpicture}
    } \subfloat[]{
      \begin{tikzpicture}[rotate=10]
	% blue filling
	\fill[draw=none,fill=blue!50!white] (0.98,-0.3482)--
        plot[smooth,tension=.5] coordinates{(.98,-.3482) (1.5,-.251)
          (2,-.06) (3,.333) (3.35,.342) (3.5,-.047) (4,-0.228)
          (5,-.348)} -- (5.98,-0.3482) -- (4.02,0.3482) --
        plot[smooth,tension=.5] coordinates{ (3.04,.3484) (2.04,.4684)
          (1.54,.6469) (1.39,1.0368) (1.04,1.0294) (0.04,.6364)
          (-.46,.4454)}-- (-0.98,0.3482) -- (0.98,-0.3482);
        \fill[draw=none,fill=blue!50!white] (0.98,-0.3482)--
        plot[smooth,tension=.5] coordinates{(.98,-.3482) (1.5,-.251)
          (2,-.06) (3,.333) (3.35,.342) } -- plot[smooth,tension=.5]
        coordinates{ (1.39,1.0368) (1.04,1.0294) (0.04,.6364)
          (-.46,.4454)}-- (-0.98,0.3482) -- (0.98,-0.3482);

        \fill[blue!50!white] (0,-1.75) arc (-90:12:-1 and 1.75) --
        (0.98,-0.3482) -- plot[smooth,tension=.5]
        coordinates{(.98,-.3482) (1.5,-.251) (2,-.06) (3,.333)
          (3.35,.342) (3.5,-.047) (4,-0.228) (5,-.348)} --
        (5.98,-0.3482) -- (5,-1.75) -- (0,-1.75); \fill[blue!50!white]
        (5,-1.75) arc (-90:-12:1 and 1.75) -- (5,0) -- (5,-1.75);
        % profile of wave
        \draw (3.35,.342) -- plot[smooth,tension=.5] coordinates{
          (1.39,1.0368) (1.04,1.0294) (0.04,.6364) (-.46,.4454)} --
        (-0.98,0.3482) -- (0.98,-0.3482)-- plot[smooth,tension=.5]
        coordinates{(.98,-.3482) (1.5,-.251) (2,-.06) (3,.333)
          (3.35,.342) (3.5,-.047) (4,-0.228) (5,-.348)} --
        (5.98,-0.3482) -- (4.02,0.3482) -- plot[smooth,tension=.5]
        coordinates{ (3.34,.3484) } ;

        % draw pipe
        \draw[->] (1.75,1.25) -- (3.25,1.25) node[midway,above]
        {$u_G$}; \draw[->] (1.75,-1.25) -- (3.25,-1.25)
        node[midway,above] {$u_L$}; \draw[-{Straight Barb[left]}]
        (3.7,-1.65) -- (4.5,-1.65) node[midway,above] {$\tau_L$};
        \draw[-{Straight Barb[right]}] (3.7,1.65) -- (4.5,1.65)
        node[midway,below] {$\tau_G$}; \draw[-{Straight Barb[left]}]
        (3.7,0.02) -- (4.5,0.02) node[midway,above=-1mm] {$\tau_I$};
        \draw[-{Straight Barb[left]}] (4.5,-0.02) -- (3.7,-0.02);
        \draw[thick] (0,0) ellipse (1 and 1.75); \draw[thick]
        (0,-1.75) -- (5,-1.75); \draw[thick] (0,1.75) -- (5,1.75);
        \draw[thick] (5,1.75) arc (90:270:-1 and 1.75);
        \draw[thick,dashed] (5,1.75) arc (90:12:-1 and 1.75);
	% \node[left=.5mm] at (0,-0.8) {$\alpha_L$};
	% \node[left=.5mm] at (0,0.8) {$\alpha_G$};

	\draw[dashed] (0,-1.75) -- (4,-2.455);
	\begin{scope}
          \path[clip] (0,-1.75) -- (4,-2.455) -- (4,-1.75);
          \draw[black] (0,-1.75) circle (35mm); \node at
          ($(0,-1.75)+(-5:25mm)$) {$\phi$};
	\end{scope}
      \end{tikzpicture}
    }
  \end{center}
  \caption{Pipe cross section and geometry}
  \label{fig: pipe cross section and geometry}
\end{figure}

We consider the four-equation two-fluid model for stratified pipe flow
with zero mass transfer. This model can be derived by averaging of the
conservation equations across the cross-sectional area (see for
example
\cite{soo1989particulates,ishii2010thermo,drew2006theory,drew1983mathematical}). The
resulting equations are commonly written
\begin{subequations}
  \begin{align*}
    \pdv{t} (\rho_K \alpha_K) + \pdv{x}(\rho_K \alpha_K v_K) &= 0,\\
    \pdv{t}(\rho_K \alpha_K v_K) + \pdv{x}(\rho_K\alpha_K v_K^2 +
    \alpha_K p_K) + \rho_K \alpha_K g_y \pdv{x}h &= s_K,
  \end{align*}
\end{subequations}
where $\rho_K$ denotes the density, $\alpha_K$ denotes the volume fraction and
$v_K$ the velocity of phase $K$ ($K$ is ``gas'' or ``liquid''). Furthermore, $h$ is the height of
the interface between the two fluids and the momentum sources are 
\begin{equation*}
  s_K= \tau_K \frac{\sigma_K}{A} \pm \tau_I \frac{\sigma_I}{A} - \rho_K
  \alpha_K g_x,
\end{equation*}
where $\tau_K$ and $\tau_I$ are the wall respectively interface shear
stress, $\sigma_K$ and $\sigma_I$ are the wetted lengths (see Figure~\ref{fig: pipe cross section and geometry}) of the $K$ phase and the interface respectively and $A$  the pipe cross-sectional area. The axial and transverse components of the gravitational acceleration are given by $g_x = g \sin\phi$ and $g_y = g\cos\phi$, where $\phi$ is the pipe inclination relative to the horizontal plane. Figure \ref{fig: pipe cross section and geometry} illustrates the pipe geometry and the aforementioned quantities.
Note that the mass $m_K = \rho_K \alpha_K$ and momentum $m_K u_K$ are conserved properties. We assume that phase $K$ is occupied by either gas, $K=G$, or liquid, $K=L$, and that the two phases are segregated from each other (see Figure \ref{fig: pipe cross section and geometry}). Note that for two phase flow
\begin{equation*}
	\alpha_G + \alpha_L = 1.
\end{equation*}
In the following we will make some simplifying assumptions. First we assume zero momentum sources and no pipe inclination. Next we assume that the gas is ideal so that $\rho_G = p_G C_G$ for some constant $C_G$. Then
\begin{equation*}
	\alpha_G p_G = \frac{m_G}{C_G}.
\end{equation*}
Lastly, assuming that the liquid is incompressible, with constant density $\rho_L$, and assuming hydrostatic balance we find
\begin{align*}
	p_L &= p_G + \frac{1}{2}\alpha_G m_G +\frac{1}{2}\alpha_L m_L\\
	& = \frac{m_G}{\alpha_G C_G} + \frac{1}{2}(1-\alpha_L)m_G + \frac{1}{2}\frac{m_L}{\rho_L} m_L\\
	& = \frac{m_G}{(1-\alpha_L)C_G} +\frac{1}{2}\left(1-\frac{m_L}{\rho_L}\right)m_G + \frac{m_L^2}{2\rho_L}\\
	& = \frac{m_G}{(1-\frac{m_L}{\rho_L})C_G} +\frac{1}{2}\left(1-\frac{m_L}{\rho_L}\right)m_G + \frac{m_L^2}{2\rho_L}.
\end{align*}
Therefore the system of equations we consider is
\begin{equation}
\begin{aligned}
	&\begin{cases}
		\pdv{t} m_G + \pdv{x}(m_G v_G) = 0,\\
		\pdv{t}(m_G v_G)+ \pdv{x}\left(m_G v_G^2 + \frac{m_G}{C_G}\right) = 0,
	\end{cases}\\
	&\begin{cases}
		\pdv{t} m_L + \pdv{x}(m_L v_L) = 0,\\
		\pdv{t} (m_L v_L) + \pdv{x}\left(m_L v_L^2 + \frac{m_L m_G}{(\rho_L -m_L)C_G} + \frac{m_L m_G}{2\rho_L^2}(\rho_L-m_L) + \frac{m_L^3}{2\rho_L^2}\right) = 0.
	\end{cases}
\end{aligned}
\label{two-phase flow problem}
\end{equation}
If we let $\bu=(m_G,m_Gv_G)$ and $\bw=(m_L,m_Lv_L)$ we see that this
can be written as
\begin{align*}
  \bu_t + f(\bu)_x &=0,\\ \bw_t + g(\bu,\bw)_x &=0,
\end{align*}
for some nonlinear functions $f$ and $g$. 
This means that the $4\times 4$ model partially decouples into two
$2\times 2$ systems, where the second system is dependent on the
first, but not vice versa. The situation is reminiscent of so-called
triangular systems of (scalar) conservation laws considered in~\cite{CMR}, in which the convergence of a finite volume scheme was
proved for such a system, in the case that $\bu$ and $\bw$ are
scalars. In our case $u$ and $w$ are vectors, but when approximating,
we can stagger the discretizations of $u$ and $w$ so that the
approximation of $u$ is continuous (constant) across each cell
interface in the approximation of $w$. The advantage of this is that
the numerical flux for $w$ can be a standard numerical flux for
$2\times 2$ system of conservation laws. In this paper we investigate
using the Roe flux as a numerical flux for $w$. While not investigated
in the present work, higher order methods can easily be built from
first order numerical fluxes by higher order reconstruction. 

The rest of this paper is organized as follows. In
Section~\ref{sec:nummeth} we detail the construction of the Roe method for the gas and the
liquid phase, as well as the Nessahu--Tadmor second order scheme. In Section~\ref{sec:numex} we describe how to calculate
the exact solution to some Riemann problems, and test our scheme on
these.

\section{Numerical Methods}\label{sec:nummeth}
In this section we will describe the numerical methods we use to find
approximate solutions to the two-phase flow problem \eqref{two-phase
  flow problem}. For the gas phase we will use Roe's method and for
the liquid phase we will use either  Roe's method or the second-order nonstaggered Nessyahu--Tadmor scheme.

\begin{figure}
\begin{center}
  \begin{tikzpicture}[scale=.5]
    \draw (0,0) -- (14,0); \foreach \x in {0,2,4,6,8,10,12,14}{ \draw
      (\x,.1) -- (\x,-.1); } \node[below] at (0,0) {$a$}; \node[below]
    at (4,0) {$x_{j-1}$}; \node[below] at (6,0) {$x_j$}; \node[below]
    at (8,0) {$x_{j+1}$}; \node[below] at (10,0) {$x_{j+2}$};
    \node[below] at (14,0) {$b$};

    \foreach \x in {2,4,6,8,10,12,14}{ \draw (\x,.1) -- (\x,-.1);
      \draw[dashed,red!40!white, thick] (\x-1,0) -- (\x-1,-5); }
    \draw[red,thick] (0,2) -- (1,2); \draw[red,thick] (1,2.5) --
    (3,2.5); \draw[red,thick] (3,3) -- (5,3); \draw[red,thick] (5,2.5)
    -- (7,2.5); \draw[red,thick] (7,2) -- (9,2); \draw[red,thick]
    (9,2.5) -- (11,2.5); \draw[red,thick] (11,2) -- (13,2);
    \draw[red,thick] (13,1.5) -- (14,1.5);
    \draw[red!40!white,thick,dashed] (1,2.5) -- (1,0);
    \draw[red!40!white,thick,dashed] (3,3) -- (3,0);
    \draw[red!40!white,thick,dashed] (5,3) -- (5,0);
    \draw[red!40!white,thick,dashed] (7,2.5) -- (7,0);
    \draw[red!40!white,thick,dashed] (9,2.5) -- (9,0);
    \draw[red!40!white,thick,dashed] (11,2.5) -- (11,0);
    \draw[red!40!white,thick,dashed] (13,2) -- (13,0); \node[above] at
    (6,2.5) {$u_j^n$}; \node[above] at (8,2) {$u_{j+1}^n$};

    \draw (0,-5) -- (14,-5); \foreach \x in {0,1,3,5,7,9,11,13,14}{
      \draw (\x,-4.9) -- (\x,-5.1); } \node[below] at (0,-5) {$a$};
    \node[below] at (5,-5) {$x_{j-\frac{1}{2}}$}; \node[below] at
    (7,-5) {$x_{j+\frac{1}{2}}$}; \node[below] at (9,-5)
    {$x_{j+\frac{3}{2}}$}; \node[below] at (14,-5) {$b$}; \draw[blue,
    thick] (0,-3.5) -- (2,-3.5); \draw[blue, thick] (2,-2.5) --
    (4,-2.5); \draw[blue, thick] (4,-3) -- (6,-3); \draw[blue, thick]
    (6,-2) -- (8,-2); \draw[blue, thick] (8,-2.5) -- (10,-2.5);
    \draw[blue, thick] (10,-3) -- (12,-3); \draw[blue, thick]
    (12,-2.5) -- (14,-2.5); \draw[blue!40!white,thick,dashed] (2,-2.5)
    -- (2,-5); \draw[blue!40!white,thick,dashed] (4,-2.5) -- (4,-5);
    \draw[blue!40!white,thick,dashed] (6,-2) -- (6,-5);
    \draw[blue!40!white,thick,dashed] (8,-2) -- (8,-5);
    \draw[blue!40!white,thick,dashed] (10,-2.5) -- (10,-5);
    \draw[blue!40!white,thick,dashed] (12,-2.5) -- (12,-5);
    \node[above] at (7,-2) {$w_{j+\frac{1}{2}}^n$};
  \end{tikzpicture}
\caption{Schematics of the numerical discretization and approximation.}\label{fig: schematics of discretization}
\end{center}
\end{figure}
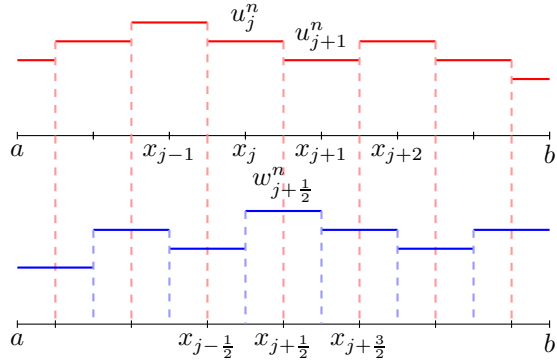

We discretize the domain $[a,b]$ using two grids that are
staggered with respect to each other. Let $\Delta x = (b-a)/N$ and $\Delta t=\frac{T}{M+1}$ and approximate $u=(m_G,m_G v_G)$ and $w=(m_L,m_L v_L)$ as follows:
\begin{align*}
  u(x_j,t^n) &\approx u_j^n & j&=0,\ldots,N,\\
  w(x_{j+\frac{1}{2}},t^n) &\approx w_{j+\frac{1}{2}}^n &
  j&=0,\ldots,N-1,
\end{align*}
where $x_j=a+j\Delta x$, $j=0,\ldots,N$ and
$x_{j+\frac{1}{2}}=a+(j+\frac{1}{2})\Delta x$, $j=0,\ldots,N-1$, as
well as $t^n=n\Delta t$, $n=0,\ldots,M$. Thus the approximation to
$u(\cdot,t^n)$ is a piecewise constant function which may have
discontinuities at $x_{j+\frac{1}{2}}$, and the approximation to $w(x,t^n)$ is a
piecewise constant function which may have discontinuities at
$x_j$. 
Figure~\ref{fig: schematics of discretization} illustrates the
staggered grid and the approximations to $u$ and $w$.

\subsection{The Roe method}\label{subsec:roe}
Roe's Riemann solver is a method to approximate the solution of the
Riemann problem for the hyperbolic system
\begin{equation*}
    \pdv{t} u + \pdv{x}f(u)  =0
\end{equation*}
with
\begin{equation*}
    u = \begin{pmatrix}
        m\\
        m v
    \end{pmatrix}\qquad\text{and}\qquad f(u) = \begin{pmatrix}
        m v\\
        m v^2 + p(m)
    \end{pmatrix}.
\end{equation*}
Here $m$ can be either $m_G$ or $m_L$ and $v$ can be $v_G$ or $v_L$
respectively, depending on which phase we consider.  Roe's
approximate Riemann solver \cite{roe1997approximate} (cf.~also
\cite{harten1983upstream} and \cite{leveque_book}) yields the solution
of the following linear problem
\begin{equation*}
    \pdv{t} u + A(u^L,u^R)\pdv{x} u = 0, \quad
    u(x,0) = \begin{cases}
        u^L & x<0,\\
        u^R & x>0.
    \end{cases}
\end{equation*}
The so-called Roe matrix $A(u^L,u^R)$ has to satisfy the following properties:
\begin{enumerate}[(1)]
    \item $A(u^L,u^R)$ is diagonizable with real eigenvalues,
    \item $A(u^L,u^R)\to f'(u)$ as $u_L,u_R\to u$ and
    \item $A(u^L,u^R)(u^R-u^L) = f(u^R)-f(u^L)$. 
\end{enumerate}
The first property is required for hyperbolicity and the second for consistency with the original conservation law. The third property ensures that single shocks of the nonlinear system are also discontinuous solutions of the linear system.
In order to find such a matrix $A=A(u^L,u^R)$ we define the parameter vector
\begin{equation*}
    z(u) = \begin{pmatrix}
        z_1\\z_2
    \end{pmatrix}
    = \begin{pmatrix}
        \sqrt{m}\\ \sqrt{m}v
    \end{pmatrix}.
\end{equation*}
Then $u$ and $f(u)$ can be written as
\begin{equation*}
    u(z) = \begin{pmatrix}
        z_1^2\\ z_1 z_2
    \end{pmatrix}
    \qquad \text{and}\qquad f(z) = \begin{pmatrix}
        z_1 z_2\\
        z_2^2 + p(z_1^2)
    \end{pmatrix}
\end{equation*}
and therefore
\begin{equation*}
    \pdv{u}{z} (z) = \begin{pmatrix}
        2 z_1 & 0\\ z_2 & z_1
    \end{pmatrix}
    \qquad \text{and}\qquad \pdv{f}{z} (z) = \begin{pmatrix}
        z_2 & z_1\\
        2z_2 p'(z_1^2) & 2 z_2
    \end{pmatrix}
\end{equation*}
Considering the straight line segment
\begin{equation*}
    \gamma(\xi) = z^L + (z^R - z^L)\xi,\qquad \xi\in[0,1],
\end{equation*}
where $z^L = z(u^L)$ and $z^R = z(u^R)$, we find
\begin{equation*}
    u^R-u^L = \int_0^1 \frac{\operatorname{d}u(\gamma(\xi))}{\dxi} \dxi
    =\int_0^1 \frac{\operatorname{d}u(\gamma(\xi))}{\operatorname{d}\gamma}z'(\xi)\dxi
    = \int_0^1  \frac{\operatorname{d}u(\gamma(\xi))}{\operatorname{d}\gamma} \dxi(z^R-z^L) =: B(z^R - z^L)
\end{equation*}
and similarly
\begin{equation*}
    f(u^R)-f(u^L) = \int_0^1 \frac{\operatorname{d}f(\gamma(\xi))}{\dxi} \dxi
    =\int_0^1 \frac{\operatorname{d}f(\gamma(\xi))}{\operatorname{d}\gamma}\gamma'(\xi)\dxi
    = \int_0^1  \frac{\operatorname{d}f(\gamma(\xi))}{\operatorname{d}\gamma} \dxi(z^R-z^L) =:  C(z^R - z^L)
\end{equation*}
Then, by setting $A = C B^{-1}$ we have
\begin{equation*}
    f(u^R) - f(u^L) = A (u^R - u^L).
\end{equation*}
Note that
\begin{equation}
    C := \int_0^1  \frac{\operatorname{d}f(z(\xi))}{\operatorname{d}\gamma} \dxi = \begin{pmatrix}
        \overline{z_2} & \overline{z_1}\\
        2\overline{z_1 p'(z_1^2)} & 2\overline{z_2}
    \end{pmatrix}
    \label{C matrix}
\end{equation}
and
\begin{equation*}
    B := \int_0^1  \frac{\operatorname{d}u(z(\xi))}{\operatorname{d}\gamma} \dxi = \begin{pmatrix}
        2 \overline{z_1} & 0\\
        \overline{z_2} & \overline{z_1}
    \end{pmatrix}
\end{equation*}
and thus
\begin{equation*}
    B^{-1} = \frac{1}{2 \overline{z_1}^2}\begin{pmatrix}
        \overline{z_1} & 0\\
        -\overline{z_2} & 2\overline{z_1}
    \end{pmatrix}
\end{equation*}
where
\begin{equation*}
    \overline{h(z_k)} = \int_0^1 h(z_k^L + (z_k^R - z_k^L)\xi) \dxi, \qquad k=1,2
\end{equation*}
for a function $h$. Therefore
\begin{equation}
    A = C B^{-1} = \frac{1}{2\overline{z_1}^2} \begin{pmatrix}
        0 & 2 \overline{z_1}^2\\
        2 \overline{z_1}\overline{z_1p'(z_1^2)} - 2\overline{z_2}^2 & 4\overline{z_1}\,\overline{z_2}
    \end{pmatrix}.
    \label{A matrix}
\end{equation}
The matrix $A$ has the eigenvalues
\begin{equation*}
    \lambda_{1,2} = \frac{\overline{z_2}}{\overline{z_1}} \mp \sqrt{\frac{\overline{z_1 p'(z_1^2)}}{\overline{z_1}}},
\end{equation*}
which are real provided that $m>0$ and $p'(m)>0$, and eigenvectors $r_{1,2} = \begin{pmatrix} 1\\ \lambda_{1,2} \end{pmatrix}$.
From this Roe's method can be defined using
\begin{equation}
    |A| = R \diag(|\lambda_1|,|\lambda_2|) R^{-1} = \frac{1}{\lambda_2-\lambda_1}\begin{pmatrix}
		\lambda_2|\lambda_1|-\lambda_1|\lambda_2| & |\lambda_2|-|\lambda_1|\\
		\lambda_2\lambda_1^2 - \lambda_1\lambda_2^2 & \lambda_2^2-\lambda_1^2
	\end{pmatrix}
	\label{A absolute}
\end{equation}
where $R=(r_1,r_2)$ is the matrix consisting of the right eigenvectors to define the numerical flux
\begin{equation}
    F^{\text{Roe}}(u^L,u^R) = \frac{1}{2}\left(f\left(u^L\right)+f\left(u^R\right)\right) - \frac{1}{2}\left|A\left(u^L,u^R\right)\right| \left(u^R-u^L\right)
    \label{Roe flux}
\end{equation}
which can then be used in the following Roe scheme:
\begin{equation*}
    u_j^{n+1} = u_j^n - \frac{\Delta t}{\Delta x}\left(F^{\text{Roe}}\left(u_j^n,u_{j+1}^n\right) - F^{\text{Roe}}\left(u_{j-1}^n,u_j^n\right)\right).
\end{equation*}

\subsubsection{The Roe method for the gas phase}
% \hili{We should really have the notation $p_{G}$ and $p_L$ here.}
The pressure in the gas phase is
\begin{equation}
	p(m_G) = \frac{m_G}{C_G} \label{pressure gas phase}
\end{equation}
with derivative
\begin{equation*}
	p'(m_G) = \frac{1}{C_G}.
\end{equation*}
Thus the Roe matrix $A$ becomes
\begin{equation*}
	A(u^L,u^R) = \begin{pmatrix}
		0 & 1\\
		-\left(\frac{\overline{z_2}}{\overline{z_1}}\right)^2 -\frac{1}{C_G} & 2\frac{\overline{z_2}}{\overline{z_1}}
	\end{pmatrix} = \begin{pmatrix}
		0 & 1\\
		-\widehat{u}^2+\frac{1}{C_G} & 2\widehat{u}
	\end{pmatrix}
\end{equation*}
where
\begin{equation*}
	\widehat{u} = \frac{\overline{z}_2}{\overline{z}_1} = \frac{\sqrt{m_G^L}v_G^L + \sqrt{m_G^R}v_G^R}{\sqrt{m_G^L}+\sqrt{m_G^R}}.
\end{equation*}
The eigenvalues of $A(u^L,u^R)$ are
\begin{equation*}
	\lambda_1 = \widehat{u}-\frac{1}{\sqrt{C_G}},\qquad\lambda_2 = \widehat{u}+\frac{1}{\sqrt{C_G}}
\end{equation*}
and hence
\begin{equation*}
	\left|A(u^L,u^R)\right| = \frac{\sqrt{C_G}}{2}\begin{pmatrix}
		\lambda_2|\lambda_1|-\lambda_1|\lambda_2| & |\lambda_2|-|\lambda_1|\\
		\lambda_2\lambda_1^2 - \lambda_1\lambda_2^2 & \lambda_2^2-\lambda_1^2
	\end{pmatrix}.
\end{equation*}

\subsubsection{The Roe method for the liquid phase}

The pressure in the liquid phase is
\begin{equation}
	P(m_G,m_L) = \frac{m_L m_G}{(\rho_L -m_L)C_G} + \frac{m_L m_G}{2\rho_L^2}(\rho_L-m_L) + \frac{m_L^3}{2\rho_L^2}
    \label{pressure liquid phase}
\end{equation}
with
\begin{equation*}
	P_{m_L}(m_G,m_L) = \frac{m_G \rho_L}{(\rho_L -m_L)^2C_G} + \frac{m_G}{2\rho_L} - \frac{m_L m_G}{\rho_L^2} + \frac{3 m_L^2}{2\rho_L^2}
\end{equation*}
In order to use the Roe scheme here, we will have to evaluate the integral
\begin{equation*}
	\overline{z_1 P_{m_L}(m_G,z_1^2)}
\end{equation*}
in \eqref{C matrix} numerically. This can be done up to machine precision with an appropriate quadrature since $P_{m_L}$ is a rational function. After calculating the eigenvalues of the resulting matrix $A(u^L,u^R)$ given by \eqref{A matrix} we can then assemble the Roe matrix $A(u^L,u^R)$ and the corresponding flux $F^{\text{Roe}}(u^L,u^R)$ as in \eqref{A absolute} respectively \eqref{Roe flux}.

\subsection{Nonstaggered second-order Nessyahu--Tadmor scheme for the
  liquid phase}
\label{subsec:ntscheme}
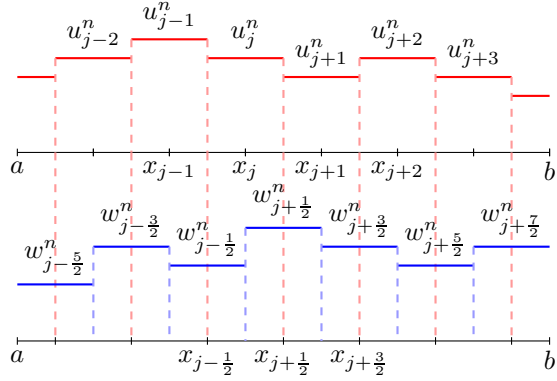
\begin{figure}
\begin{center}
	\begin{tikzpicture}[scale=.5]
		\draw (0,0) -- (14,0);
		\foreach \x in {0,2,4,6,8,10,12,14}{
		\draw (\x,.1) -- (\x,-.1);
		}
		\node[below] at (0,0) {$a$};
		\node[below] at (4,0) {$x_{j-1}$};
		\node[below] at (6,0) {$x_j$};
		\node[below] at (8,0) {$x_{j+1}$};
		\node[below] at (10,0) {$x_{j+2}$};
		\node[below] at (14,0) {$b$};

		\foreach \x in {2,4,6,8,10,12,14}{
		\draw (\x,.1) -- (\x,-.1);
		\draw[dashed,red!40!white,thick] (\x-1,0) -- (\x-1,-5);
		}
		\draw[red, thick] (0,2) -- (1,2);
		\draw[red, thick] (1,2.5) -- (3,2.5);
		\draw[red, thick] (3,3) -- (5,3);
		\draw[red, thick] (5,2.5) -- (7,2.5);
		\draw[red, thick] (7,2) -- (9,2);
		\draw[red, thick] (9,2.5) -- (11,2.5);
		\draw[red, thick] (11,2) -- (13,2);
		\draw[red, thick] (13,1.5) -- (14,1.5);
		\draw[red!40!white,thick,dashed] (1,2.5) -- (1,0);
		\draw[red!40!white,thick,dashed] (3,3) -- (3,0);
		\draw[red!40!white,thick,dashed] (5,3) -- (5,0);
		\draw[red!40!white,thick,dashed] (7,2.5) -- (7,0);
		\draw[red!40!white,thick,dashed] (9,2.5) -- (9,0);
		\draw[red!40!white,thick,dashed] (11,2.5) -- (11,0);
		\draw[red!40!white,thick,dashed] (13,2) -- (13,0);
		\node[above] at (2,2.5) {$u_{j-2}^n$};
		\node[above] at (4,3) {$u_{j-1}^n$};
		\node[above] at (6,2.5) {$u_j^n$};
		\node[above] at (8,2) {$u_{j+1}^n$};
		\node[above] at (10,2.5) {$u_{j+2}^n$};
		\node[above] at (12,2) {$u_{j+3}^n$};

		\draw (0,-5) -- (14,-5);
		\foreach \x in {0,1,3,5,7,9,11,13,14}{
		\draw (\x,-4.9) -- (\x,-5.1);
		}
		\node[below] at (0,-5) {$a$};
		\node[below] at (5,-5) {$x_{j-\frac{1}{2}}$};
		\node[below] at (7,-5) {$x_{j+\frac{1}{2}}$};
		\node[below] at (9,-5) {$x_{j+\frac{3}{2}}$};
		\node[below] at (14,-5) {$b$};
		\draw[blue, thick] (0,-3.5) -- (2,-3.5);
		\draw[blue, thick] (2,-2.5) -- (4,-2.5);
		\draw[blue, thick] (4,-3) -- (6,-3);
		\draw[blue, thick] (6,-2) -- (8,-2);
		\draw[blue, thick] (8,-2.5) -- (10,-2.5);
		\draw[blue, thick] (10,-3) -- (12,-3);
		\draw[blue, thick] (12,-2.5) -- (14,-2.5);
		\draw[blue!40!white,thick,dashed] (2,-2.5) -- (2,-5);
		\draw[blue!40!white,thick,dashed] (4,-2.5) -- (4,-5);
		\draw[blue!40!white,thick,dashed] (6,-2) -- (6,-5);
		\draw[blue!40!white,thick,dashed] (8,-2) -- (8,-5);
		\draw[blue!40!white,thick,dashed] (10,-2.5) -- (10,-5);
		\draw[blue!40!white,thick,dashed] (12,-2.5) -- (12,-5);
		\node[above] at (1,-3.5) {$w_{j-\frac{5}{2}}^n$};
		\node[above] at (3,-2.5) {$w_{j-\frac{3}{2}}^n$};
		\node[above] at (5,-3) {$w_{j-\frac{1}{2}}^n$};
		\node[above] at (7,-2) {$w_{j+\frac{1}{2}}^n$};
		\node[above] at (9,-2.5) {$w_{j+\frac{3}{2}}^n$};
		\node[above] at (11,-3) {$w_{j+\frac{5}{2}}^n$};
		\node[above] at (13,-2.5) {$w_{j+\frac{7}{2}}^n$};
	\end{tikzpicture}
\caption{Stencil of the numerical method for the liquid phase.}\label{fig: stencil}
\end{center}
\end{figure}

We will now discuss a second-order scheme for the liquid phase. We want to solve the conservation law
\begin{equation*}
    \frac{\partial}{\partial t} w + \frac{\partial}{\partial x} g(w,u) =0
\end{equation*}
with
\begin{equation*}
    w=\begin{pmatrix}
        m_L\\ m_L v_L
    \end{pmatrix},
\qquad
g(w,u) = \begin{pmatrix}
    m_L v_L\\ m_L v_L^2 + \frac{m_Lm_G}{(\rho_L-m_L)C_G} + \frac{m_Lm_G}{2\rho_L^2}(\rho_L-m_L) + \frac{m_L^3}{2\rho_L^2}
\end{pmatrix}
\end{equation*}
and $u$ as before.
Starting from the (staggered) second-order Nessyahu--Tadmor scheme \cite{nessyahu1990non}
\begin{align*}
	w_{j+\frac{1}{2}}^{n+\frac{1}{2}} &= w_{j+\frac{1}{2}}^n - \frac{\lambda}{2}g'_{j+\frac{1}{2}},\\
	w_j^{n+1} &= \frac{1}{2}\left(w_{j+\frac{1}{2}}^n+w_{j-\frac{1}{2}}^n\right) + \frac{1}{8}\left(w'_{j-\frac{1}{2}}-w'_{j+\frac{1}{2}}\right) -\lambda\left(g\left(w_{j+\frac{1}{2}}^{n+\frac{1}{2}}\right)-g\left(w_{j-\frac{1}{2}}^{n+\frac{1}{2}}\right)\right).
\end{align*}
where $\lambda=\frac{\Delta t}{\Delta x}$ and the discrete derivatives $g_{j+\frac{1}{2}}'$ and $w_{j+\frac{1}{2}}'$ are made precise below.
We use the averaging procedure described in \cite{jiang1998high} to get a nonstaggered version. To this end we reconstruct a piecewise-linear interpolant through the staggered cell-averages at time $t^{n+1}$:
\begin{equation*}
	L_{j}^{n+1}(x)= w_j^{n+1} + w'_{j}\left(\frac{x-x_j}{\Delta x}\right),\qquad x\in\left(x_{j-\frac{1}{2}},x_{j+\frac{1}{2}}\right)
\end{equation*}
where the staggered discrete derivative $w'_j$ is given by
\begin{equation}
	w'_j = \minmod\left(\Delta w_{j+\frac{1}{2}}^{n+1},\Delta w_{j-\frac{1}{2}}^{n+1}\right) \label{vprimej}
\end{equation}
with
\begin{align*}
	\Delta w_{j+\frac{1}{2}}^{n+1} &= w_{j+1}^{n+1} - w_j^{n+1}\\
	&=\frac{1}{2}\left( w_{j+\frac{3}{2}}^n - w_{j-\frac{1}{2}}^n\right) -\frac{1}{8}\left( w'_{j+\frac{3}{2}}-2w'_{j+\frac{1}{2}}+w'_{j-\frac{1}{2}} \right) - \lambda \left( g\left(w_{j+\frac{3}{2}}^{n+\frac{1}{2}}\right) -2g\left(w_{j+\frac{1}{2}}^{n+\frac{1}{2}} \right) +g\left(w_{j-\frac{1}{2}}^{n+\frac{1}{2}} \right) \right)
\end{align*}
Here, the minmod limiter is defined as
\begin{equation*}
  \minmod(a_1,\ldots,a_n) = \begin{cases}
    \sign(a_1)\min_{1\leq k\leq n}|a_k| &\text{if }\sign(a_1)=\ldots=\sign(a_n),\\
    0 &\text{otherwise.}
  \end{cases}
\end{equation*}
We then average these interpolants over the cell $(x_j,x_{j+1})$ to obtain a nonstaggered scheme:
\begin{align*}
	w_{j+\frac{1}{2}}^{n+1} &= \frac{1}{\Delta x}\left( \int_{x_j}^{x_{j+\frac{1}{2}}}L_j^{n+1}(x)  dx +\int_{x_{j+\frac{1}{2}}}^{x_{j+1}} L_{j+1}^{n+1}(x) dx \right)\\
	&= \frac{1}{2}\left( w_j^{n+1}+w_{j+1}^{n+1} \right) -\frac{1}{8}\left( w'_{j+1} -w'_j \right)\\
	&= \frac{1}{4}\left(w_{j+\frac{3}{2}}^n+2w_{j+\frac{1}{2}}^n + w_{j-\frac{1}{2}}^n\right) - \frac{1}{16}\left(w'_{j+\frac{3}{2}}-w'_{j-\frac{1}{2}}\right) \\
	&\phantom{=} -\frac{\lambda}{2}\left( g\left(w_{j+\frac{3}{2}}^{n+\frac{1}{2}} \right) - g\left(w_{j-\frac{1}{2}}^{n+\frac{1}{2}} \right) \right) - \frac{1}{8}\left( w'_{j+1}-w'_j \right).
\end{align*}
Here, $(w'_j)$ and $(w'_{j+\frac{1}{2}})$ are given by \eqref{vprimej} and 
\begin{equation*}
	w'_{j+\frac{1}{2}} = \minmod\left( \Delta w_{j+1}^n,\frac{1}{2}\left( \Delta w_{j+1}^n+ \Delta w_j^n \right),\Delta w_j^n \right)
\end{equation*}
respectively, where $\Delta w_j^n = w_{j+\frac{1}{2}}^n-w_{j-\frac{1}{2}}^n$. Note that, in our application, the flux $g$ depends on both $w$ and $u$, i.e.,
\begin{equation*}
	g\left(w_{j+\frac{1}{2}}^{n+\frac{1}{2}}\right) = g\left( w_{j+\frac{1}{2}}^n -\frac{\lambda}{2}g'_{j+\frac{1}{2}},u_{j+\frac{1}{2}}^n\right),
\end{equation*}
where the discrete derivatives of the flux are given by
\begin{equation*}
	g'_{j+\frac{1}{2}} = \minmod\left( g\left(w_{j+\frac{3}{2}}^n,u_{j+1}^n\right)-g\left( w_{j+\frac{1}{2}}^n,u_{j+1}^n \right),g\left( w_{j+\frac{1}{2}}^n,u_j^n \right) - g\left( w_{j-\frac{1}{2}}^n,u_j^n \right) \right)
\end{equation*}
and $u_{j+\frac{1}{2}}^n$ denotes the solution of the Riemann problem
\begin{equation*}
  \begin{gathered}
    \frac{\partial}{\partial t} u + A(u_j^n,u_{j+1}^n)\frac{\partial}{\partial x} u =0,\qquad  t>t^n,\\
    u(x,t^n) = \begin{cases}
      u_j^n,&x<x_{j+\frac{1}{2}}\\
      u_{j+1}^n, & x>x_{j+\frac{1}{2}}
    \end{cases}
  \end{gathered}
\end{equation*}
discussed in the previous section evaluated at $x=x_{j+\frac{1}{2}}$, i.e.,
\begin{equation*}
	u_{j+\frac{1}{2}}^n =\begin{cases}
		u_j^n, & \text{if } \lambda_1,\lambda_2>0\\
		\frac{\sqrt{C_G}}{2}\left( (\lambda_2 u_{j+1}^1-u_{j+1}^2)r_1 + (-\lambda_1 u_j^1+u_j^2)r_2 \right), & \text{if } \lambda_1<0<\lambda_2\\
		u_{j+1}^n, &\text{if } \lambda_1,\lambda_2<0,
	\end{cases}
\end{equation*}
%%%% See script from Sid chapter 6.3
where $u_k^1$ and $u_k^2$ denote the first and second component of $u_k$ respectively (see \cite{leveque_book}).
Figure \ref{fig: stencil} illustrates all the values of $(u_j^n)$ and $(w_{j+\frac{1}{2}}^n)$ that are used when calculating $w_{j+\frac{1}{2}}^{n+1}$.

Note that the same averaging procedure detailed above can also be applied for higher-order schemes, cf.~\cite{jiang1998high}.

\section{Numerical experiments}\label{sec:numex}
In order to have available exact solutions to which we can compare our
approximations, we consider Riemann problems, i.e., the initial value
problem for \eqref{two-phase flow problem} where the initial data
consists of a single jump between two constant values, viz.,
\begin{equation*}
  \begin{pmatrix}
    m_G\\
    v_G\\
    m_L\\
    v_L
  \end{pmatrix}(0,x) = \begin{cases}
    \begin{pmatrix}
      m_G^L\\
      v_G^L\\
      m_L^L\\
      v_L^L
    \end{pmatrix},& \text{if }x<0,\\
    \begin{pmatrix}
      m_G^R\\
      v_G^R\\
      m_L^R\\
      v_L^R
    \end{pmatrix},& \text{if }x>0.
  \end{cases}
\end{equation*}
If we fix the left state, we can find Riemann problems and their
solution by following solution curves (rarefaction curves and Hugoniot
loci) in phase space, see \cite{HR_book}. In the following sections we
will first detail the ingredients necessary to find  Riemann problems
and their entropy solutions, and then present some example test cases and finally compare our schemes to those test cases.

\subsection{Ingredients to solve the Riemann problem}\label{Finding an exact solution}
We will now describe how one can find solutions to certain Riemann
problems in two cases.
Define
\begin{align*}
  p(m_G)&=\frac{m_G}{C_G}, \quad \text{and}\\
  P(m_G,m_L)&=\frac{m_Lm_G}{(\rho_L-m_L)C_G} +
  \frac{m_Lm_G}{2\rho_L^2}(\rho_L-m_L) + \frac{m_L^3}{2\rho_L^2}.
\end{align*}
and introduce the variables $q_L = m_Lv_L$ and $q_G = m_Gv_G$, and the
flux function
\begin{equation*}
  F
  \begin{pmatrix}
    m_G \\ q_g \\ m_L\\ q_L
  \end{pmatrix}
= \begin{pmatrix}
    q_G\\
    \frac{q_G^2}{m_G} + p(m_G)\\
    q_L\\
    \frac{q_L^2}{m_L} + P(m_G,m_L)
  \end{pmatrix}
\end{equation*}
so that \eqref{two-phase flow problem} reads
\begin{equation*}
  \frac{\partial}{\partial t} 
  \begin{pmatrix}
    m_G \\ q_G \\ m_L\\ q_L
  \end{pmatrix} + \frac{\partial}{\partial x}
   F
  \begin{pmatrix}
    m_G\\ q_g\\ m_L\\ q_L
  \end{pmatrix} = 0.
\end{equation*}
The Jabobian of $F$ reads
\begin{equation}
  \begin{pmatrix}
    0 & 1 & 0 & 0\\
    p'(m_G)-v_G^2 & 2v_G & 0 & 0\\
    0 & 0 & 0 & 1\\
    P_{m_G}(m_G,m_L) & 0 & P_{m_L}(m_G,m_L) - v_L^2 & 2v_L
  \end{pmatrix},
  \label{Jacobian}
\end{equation}
where $P_{m_L}=\frac{\partial P}{\partial m_L}$ and
$P_{m_G}=\frac{\partial P}{\partial m_G}$.
This matrix has  eigenvalues
\begin{align*}
  \lambda_1 &= v_G - \sqrt{p'(m_G)},& \lambda_2 &= v_G + \sqrt{p'(m_G)},\\
  \mu_1 &= v_L - \sqrt{P_{m_L}(m_G,m_L)}, & \mu_2 &= v_L +
  \sqrt{P_{m_L}(m_G,m_L)}.
\end{align*}
Observe that this system is not strictly hyperbolic since the
eigenvalues of the $\mu$ families can coincide with the eigenvalues of
the $\lambda$ families.

\subsubsection{Shocks}
In the following, we will let
\begin{equation*}
	\llbracket a \rrbracket = a^R - a^L
\end{equation*}
denote the jump in some quantity $a$.
Given $m^L,m^R$ and either $v^L$ or $v^R$ the Rankine-Hugoniot loci in either phase can be computed by solving
\begin{align}
	s \llbracket m \rrbracket &= \llbracket m v \rrbracket\label{Rankine-Hugoniot condition I}\\
	s \llbracket m v \rrbracket &= \llbracket m v^2 \rrbracket + \llbracket p \rrbracket\label{Rankine-Hugoniot condition II}
\end{align}
for $v^R$ respectively $v^L$, i.e. solving
\begin{equation}
	\llbracket v \rrbracket^2 = \frac{1}{m^R m^L} \llbracket m\rrbracket \llbracket p\rrbracket. \label{Rankine-Hugoniot equation}
\end{equation}
Here, the shock speed is
\begin{equation*}
	s = \frac{\llbracket m v\rrbracket}{\llbracket m\rrbracket}.
\end{equation*}

Note that all four variables change over shocks of the $\lambda$ families. More specifically, rearranging the terms in \eqref{Rankine-Hugoniot condition II} for the liquid phase we find
\begin{align}
    s \llbracket m_L \rrbracket &= \llbracket m_L v_L \rrbracket & &\Leftrightarrow &  m_L^R (s-v_L^R) &= m_L^L(s-v_L^L)\label{first equivalence}
    \intertext{and also}
    s \llbracket m_L \rrbracket &= \llbracket m_L v_L \rrbracket & &\Leftrightarrow & \llbracket v_L\rrbracket &= \frac{1}{m_L^R} (s-v_L^L)\llbracket m_L\rrbracket
    \label{second equivalence}
\end{align}
Hence, using \eqref{first equivalence} and \eqref{second equivalence} in \eqref{Rankine-Hugoniot condition II} we get
\begin{align*}
    \llbracket P(\cdot,\cdot) \rrbracket &= s\llbracket m_L v_L\rrbracket - \llbracket m_L v_L^2\rrbracket\\
    &= v_L^R m_L^R (s-v_L^R) - v_L^L m_L^L(s-v_L^L)\\
    &= m_L^L(s-v_L^L)\llbracket v_L\rrbracket\\
    &= \frac{m_L^L}{m_L^R} (s-v_L^L)^2\llbracket m_L\rrbracket
\end{align*}
Here, $\llbracket P(\cdot,\cdot)\rrbracket = P(m_G^R,m_L^R) - P(m_G^L,m_L^L)$.

Thus given a shock in the gas phase separating $(m_G^L,v_G^L)$ and $(m_G^R,v_G^R)$ traveling with speed $s$, and values $m_L^L$ and $v_L^L$ we can find the right state $(m_L^R,v_L^R)$ by solving
\begin{equation*}
	H(m_L^R; m_L^L, m_G^L, m_G^R, v_L^L, s) := \frac{m_L^R}{m_L^L}\llbracket P(\cdot,\cdot) \rrbracket - \llbracket m_L \rrbracket (s-v_L^L)^2 = 0
\end{equation*}
for $m_L^R$ and computing the corresponding $v_L^R$ with \eqref{second equivalence}.

\subsubsection{Rarefactions of the $\mu$ families}
Since the eigenvectors of the Jacobian \eqref{Jacobian} associated with $\mu_{1,2}$ are
\begin{equation*}
    r_{1,2}^\mu = \begin{pmatrix}
        0\\0\\ 1\\ \mu_{1,2}
    \end{pmatrix}
\end{equation*}
Rarefaction curves in the liquid phase are solutions of
\begin{equation*}
    \dv{}{\xi} \begin{pmatrix} m_L\\q_L\end{pmatrix} = \begin{pmatrix}
        1\\ \mu_{1,2}
    \end{pmatrix}
\end{equation*}
or
\begin{equation*}
    \dv{q_L}{m_L} = \mu_{1,2} = \frac{q_L}{m_L} \mp \sqrt{P_{m_L}(m_G,m_L)}
\end{equation*}
If we use $v_L = v_L(m_L,q_L)$ instead, we get
\begin{equation*}
	\dv{v_L}{m_L} = \pdv{v_L}{m_L} + \pdv{v_L}{q_L}\dv{q_L}{m_L} = \pm \frac{\sqrt{P_{m_L}(m_G,m_L)}}{m_L}
\end{equation*}
From this we see that rarefaction curves through a point $(m_L^L,q_L^L)$ are given by
\begin{equation}
	v_L^R = v_L^L \pm \int_{m_L^L}^{m_L^R} \frac{\sqrt{P_{m_L}(m_G,\xi)}}{\xi} d\xi. \label{Equation for v_L^R rarefaction}
\end{equation}

\subsubsection{Rarefactions of the $\lambda$ families}
The eigenvectors associated with $\lambda_{1,2}$, after renormalizing such that $\nabla \lambda_{1,2}\cdot r_{1,2}^\lambda = 1$, are
\begin{equation*}
  r_{1,2}^\lambda = \pm \left( \frac{p''(m_G)}{2\sqrt{p'(m_G)}} + \frac{\sqrt{p'(m_G)}}{m_G} \right) \begin{pmatrix}
    1\\ \lambda_{1,2}\\ c(m_G,v_G,m_L,v_L) \\ c(m_G,v_G,m_L,v_L)\lambda_{1,2}
  \end{pmatrix}
\end{equation*}
where
\begin{equation*}
    c(m_G,v_G,m_L,v_L) = \frac{P_{m_G}(mg,m_L)}{\lambda_{1,2}^2-P_{m_L}(m_G,m_L)+v_L^2-2v_L\lambda_{1,2}}.
\end{equation*}
By solving
\begin{equation}
    \dv{}{\xi} \begin{pmatrix}
        m_G\\ q_G\\ m_L\\ q_L
    \end{pmatrix}(\xi)
    = r_{1,2}^\lambda
    \label{big rarefaction}
\end{equation}
with the initial condition
\begin{equation*}
    \begin{pmatrix}
        m_G\\ q_G\\ m_L\\ q_L
    \end{pmatrix}(\lambda_{1,2}(m_G^L,v_G^L)) = \begin{pmatrix}
        m_G^L\\ q_G^L\\ m_L^L\\ q_L^L
    \end{pmatrix}
\end{equation*}
we get a rarefaction wave in all four variables.

\subsection{Finding suitable Riemann problems and their solutions}

We will now formulate two algorithms to find a Riemann problem whose
solution consists of a shock in the liquid phase, followed by a shock
over which all four components change, followed by a second shock in
the liquid phase, and a Riemann problem whose solution consists of a
rarefaction wave in the liquid phase, followed by a rarefaction wave
over which all four components change, followed by a second rarefaction
wave in the liquid phase respectively.
\paragraph{Algorithm 1: All-shock Riemann solution}
\begin{enumerate}[(1)]
	\item Choose parameters $\begin{pmatrix}m_G^L\\ v_G^L\end{pmatrix}, m_G^R$ as well as $m_L^L,m_L'$ and $m_L^R$.
	\item Find $\textcolor{blue}{v_G^R}$ by solving the Rankine-Hugoniot equation \eqref{Rankine-Hugoniot equation} in the liquid phase such that the shock separating the left state $\begin{pmatrix}m_G^L\\ v_G^L\end{pmatrix}$ and right state $\begin{pmatrix}m_G^R\\ \textcolor{blue}{v_G^R}\end{pmatrix}$ (with speed $s$) is a Lax shock of the $\lambda_1$ family.
	\item Find $\textcolor{green!60!black}{v_L^L}$ by solving the Rankine-Hugoniot equation \eqref{Rankine-Hugoniot equation} in the liquid phase such that the shock separating the left state $\begin{pmatrix}m_L^L\\ \textcolor{green!60!black}{v_L^L}\end{pmatrix}$ and right state $\begin{pmatrix}m_L'\\ v_L'\end{pmatrix}$ (with speed $s_L^L$ such that $s_L^L<s$) is a Lax shock of the $\mu_1$ family. Here, $\llbracket P \rrbracket = p(m_G^L,m_L') - p(m_G^L,m_L^L)$.\label{Find left-going shock}
	\item Solve $H(\textcolor{red}{m_L''}, m_L', m_G^L, m_G^R,
          v_L',s)=0$ for $\textcolor{red}{m_L''}$ and find the
          corresponding $v_L''$ with equation \eqref{second equivalence} so that the liquid phase has a shock with the same speed as the one in the gas phase.
	\item Check that $s<\mu_2(m_L'',v_L'',m_G^R)$. If not, choose different values for the parameters and start again.
	\item Find $\textcolor{teal}{v_L^R}$ by solving the Rankine-Hugoniot equation \eqref{Rankine-Hugoniot equation} such that the shock separating the left state $\begin{pmatrix}m_L''\\ v_L''\end{pmatrix}$ and right state $\begin{pmatrix}m_L^R\\ \textcolor{teal}{v_L^R}\end{pmatrix}$ (with speed $s_L^R > s$) is a Lax shock of the $\mu_2$ family. Here, $\llbracket P \rrbracket = p(m_G^R,m_L^R) - p(m_G^R,m_L'')$.
	% \item Check that $\mu_2(m_L^R,v_L^R,m_G^R) > \mu_2(m_L'',v_L'',m_G^R)$. If not, choose different values for the parameters and start again.
	% \item Find $\textcolor{teal}{v_L^R}$ by solving the Rankine-Hugoniot equation \eqref{Rankine-Hugoniot equation} in the liquid phase such that the shock separating the left state $\begin{pmatrix}m_L''\\ v_L''\end{pmatrix}$ and right state $\begin{pmatrix}m_L^R\\ \textcolor{teal}{v_L^R}\end{pmatrix}$ (with speed $s_L^R$ such that $s_L^R > s$) is a Lax shock of the $\mu_+$ family. Here, $\llbracket p \rrbracket = p(m_G^R,m_L') - p(m_G^R,m_L^L)$.
\end{enumerate}
Following this procedure we arrive at the following constant states in an all-shock solution of the Riemann problem:
\begin{equation*}
	\begin{pmatrix}
		m_G^L\\
		v_G^L\\
		m_L^L\\
		v_L^L
	\end{pmatrix}
	\underset{s_L^L}{\overset{\mu_1}{\longrightarrow}}
	\begin{pmatrix}
		m_G^L\\
		v_G^L\\
		m_L'\\
		v_L'
	\end{pmatrix}
	\underset{s}{\overset{\lambda_1}{\longrightarrow}}
	\begin{pmatrix}
		m_G^R\\
		v_G^R\\
		m_L''\\
		v_L''
	\end{pmatrix}
	\underset{s_L^R}{\overset{\mu_2}{\longrightarrow}}
	\begin{pmatrix}
		m_G^R\\
		v_G^R\\
		m_L^R\\
		v_L^R
	\end{pmatrix}
\end{equation*}
It is straightforward to extend this algorithm to include two shock waves associated with $\lambda_1$ and $\lambda_2$ by repeating Algorithm 1 for a shock of the $\lambda_2$ family and paying close attention to the ordering of the wave speeds.

\paragraph{Algorithm 2: All-rarefaction Riemann solution:}
\begin{enumerate}[(1)]
    \item Chose parameters $\begin{pmatrix}m_G^L\\v_G^L\end{pmatrix},\begin{pmatrix}m_L'\\v_L'\end{pmatrix}$ and $m_L^L$ as well as $m_L^R$.
    \item Find $\begin{pmatrix}m_G^R\\v_G^R\end{pmatrix}$ and $\begin{pmatrix}m_L''\\v_L''\end{pmatrix}$ by solving \eqref{big rarefaction} to get a rarefaction over which all four components change of the $\lambda_1$ family.
    \item Find $\textcolor{green!60!black}{v_L^L}$ by solving \eqref{Equation for v_L^R rarefaction} such that we get the rarefaction wave connecting the left state $\begin{pmatrix}m_L^L\\ \textcolor{green!60!black}{v_L^L}\end{pmatrix}$ and right state $\begin{pmatrix}m_L'\\v_L'\end{pmatrix}$ is of the $\mu_1$ family.
    \item Find $\textcolor{teal}{v_L^R}$ by solving \eqref{Equation for v_L^R rarefaction} such that we get the rarefaction wave connecting the left state $\begin{pmatrix}m_L''\\ v_L''\end{pmatrix}$ and right state $\begin{pmatrix}m_L^R\\\textcolor{teal}{v_L^R}\end{pmatrix}$ is of the $\mu_2$ family.
\end{enumerate}
Following this procedure we arrive at the following constant states in an all-rarefaction solution of the Riemann problem:
\begin{equation*}
    \begin{pmatrix}
        m_G^L\\
        v_G^L\\
        m_L^L\\
        v_L^L
    \end{pmatrix}
    \overset{\mu_1}{\longrightarrow}
    \begin{pmatrix}
        m_G^L\\
        v_G^L\\
        m_L'\\
        v_L'
    \end{pmatrix}
    \overset{\lambda_1}{\longrightarrow}
    \begin{pmatrix}
        m_G^R\\
        v_G^R\\
        m_L''\\
        v_L''
    \end{pmatrix}
    \overset{\mu_2}{\longrightarrow}
    \begin{pmatrix}
        m_G^R\\
        v_G^R\\
        m_L^R\\
        v_L^R
    \end{pmatrix}
\end{equation*}
It is straightforward to extend this algorithm to include two rarefaction waves associated with $\lambda_1$ and $\lambda_2$ or combine it with the algorithm above to include a $\lambda_2$ shock.

\paragraph{Experiment 1: All-shock Riemann solution}
By following Algorithm 1 we can create our first test case:
\begin{align*}
    \begin{pmatrix}
        m_G^L\\
        v_G^L
    \end{pmatrix}
    = \begin{pmatrix}
        2\\
        1.5
    \end{pmatrix}
    &\underset{s}{\overset{\lambda_1}{\longrightarrow}}
    \begin{pmatrix}
        m_G^R\\
        v_G^R
    \end{pmatrix}
    = \begin{pmatrix}
        2.5\\
        1.2764
    \end{pmatrix}
\\
\begin{pmatrix}
        m_L^L\\
        v_L^L
    \end{pmatrix}
    = \begin{pmatrix}
        3\\
        1
    \end{pmatrix}
    \underset{s_L^L}{\overset{\mu_1}{\longrightarrow}}
    \begin{pmatrix}
        m_L'\\
        v_L'
    \end{pmatrix}
    = \begin{pmatrix}
        3.25\\
        0.7487
    \end{pmatrix}
    &\underset{s}{\overset{\lambda_1}{\longrightarrow}}
    \begin{pmatrix}
        m_L''\\
        v_L''
    \end{pmatrix}
    = \begin{pmatrix}
        3.4995\\
        0.7226
    \end{pmatrix}
    \underset{s_L^R}{\overset{\mu_2}{\longrightarrow}}
    \begin{pmatrix}
        m_L^R\\
        v_L^R
    \end{pmatrix}
    = \begin{pmatrix}
        3\\
        0.2475
    \end{pmatrix}.
\end{align*}
Here we have used $p(m_G)$ and $p(m_G,m_L)$ as defined in \eqref{pressure gas phase} respectively \eqref{pressure liquid phase} with $C_G = 1$ and $\rho_L = 1$.
The solution to this Riemann problem consists of a shock in the liquid phase followed by a shock over which all four components change followed by a second shock in the liquid phase.
Specifically, we have
\begin{gather*}
    \mu_1(m_L^L,v_L^L,m_G^L) = -2 > s_L^L \approx -2.2667 > \mu_1(m_L',v_L',m_G^L) \approx -2.5283,\\
    \lambda_1(m_G^L,v_G^L) = 0.5 > s \approx 0.3820 > \lambda_1(m_G^R,v_G^R) \approx 0.2764,\\
    \mu_2(m_L'',v_L'',m_G^R) \approx 4.0798 > s_L^R \approx 3.5761 >  \mu_2(m_L^R,v_L^R,m_G^R) \approx 3.0537\\
\end{gather*}
so all shocks are Lax shocks.

\paragraph{Experiment 2: All-rarefaction Riemann solution}
Following Algorithm 2 we arrive at the second test case:
\begin{align*}
    \begin{pmatrix}
        m_G^L\\
        v_G^L\\
    \end{pmatrix}
    = \begin{pmatrix}
        0.4\\
        1.5\\
    \end{pmatrix}
    &\stackrel{\lambda_1}{\longrightarrow}
    \begin{pmatrix}
        m_G^R\\
        v_G^R\\
    \end{pmatrix}
    = \begin{pmatrix}
        0.2963\\
        1.8\\
    \end{pmatrix}
\\
\begin{pmatrix}
        m_L^L\\
        v_L^L
    \end{pmatrix}
    = \begin{pmatrix}
        0.7\\
        0.4141
    \end{pmatrix}
    \stackrel{\mu_1}{\longrightarrow}
    \begin{pmatrix}
        m_L'\\
        v_L'
    \end{pmatrix}
    = \begin{pmatrix}
        0.5\\
        1
    \end{pmatrix}
    &\stackrel{\lambda_1}{\longrightarrow}
    \begin{pmatrix}
        m_L''\\
        v_L''
    \end{pmatrix}
    = \begin{pmatrix}
        0.5695\\
        0.9566
    \end{pmatrix}
    \stackrel{\mu_2}{\longrightarrow}
    \begin{pmatrix}
        m_L^R\\
        v_L^R
    \end{pmatrix}
    = \begin{pmatrix}
        0.7\\
        1.3021
    \end{pmatrix}.
\end{align*}
Again we used $p(m_G)$ and $P(m_G,m_L)$ as defined in \eqref{pressure gas phase} respectively \eqref{pressure liquid phase} with $C_G = 1$ and $\rho_L = 1$. To solve the ODE \eqref{big rarefaction} numerically we used a very high-order Runge--Kutta method and evaluated the integrals \eqref{Equation for v_L^R rarefaction} using appropriate quadrature rules.
Since
\begin{gather*}
    \mu_1(m_L^L,v_L^L,m_G^L) \approx -1.8441 < \mu_1(m_L',v_L',m_G^L)\approx -0.4053 \\
    \lambda_1(m_G^L,v_G^L) = 0.5 < \lambda_1(m_G^R,v_G^R) = 0.8\\
    \mu_2(m_L'',v_L'',m_G^R) \approx 2.3936 < \mu_2(m_L^R,v_L^R,m_G^R) \approx 3.2941
\end{gather*}
the rarefaction waves are ordered as described above.

\subsection{Comparison with numerical solutions}
Figures~\ref{fig: Riemann problem 1} and \ref{fig: Riemann problem 2} show the exact and numerical solutions of Experiment 1 and 2 respectively with open boundaries.
Here we divided the spatial domain $[-5,5]$ into $54$ respectively
$55$ ($N=55$) cells and use the time discretization parameter $\Delta t=\frac{1}{150}$
to approximate the solution at $T=1$.
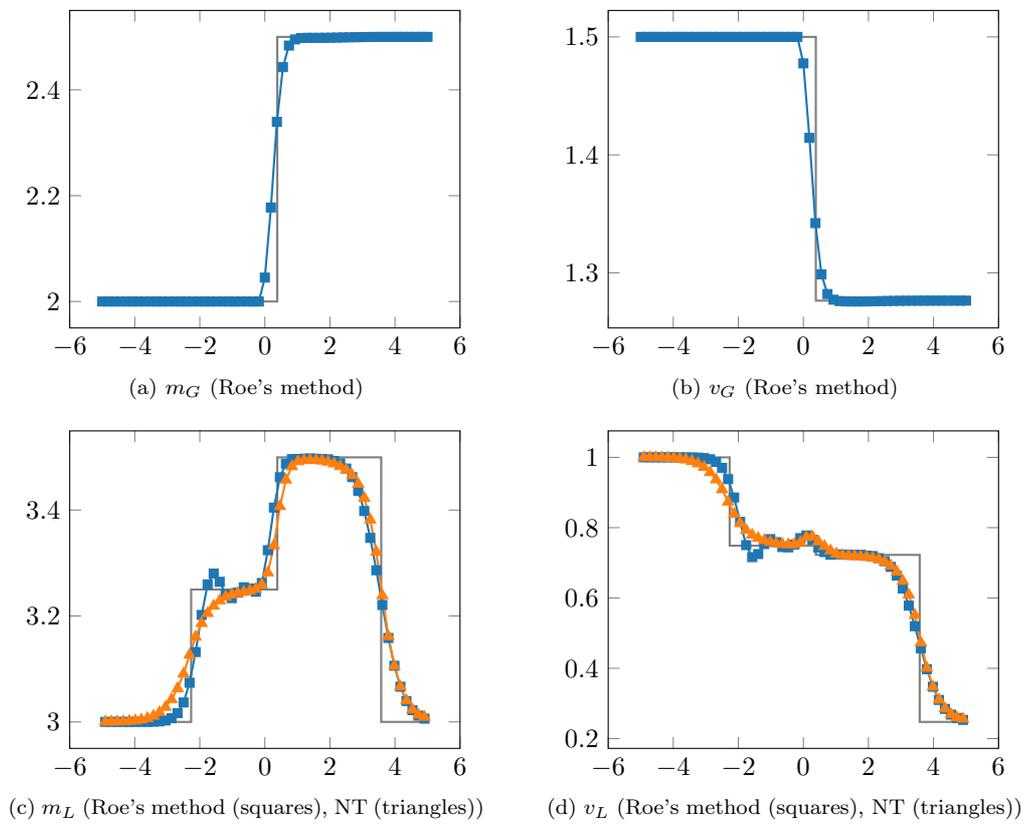
\begin{figure}
	\centering
	\subfloat[$m_G$ (Roe's method)]{
	\begin{tikzpicture}
		\begin{axis}
		\addplot[gray, thick, mark=none] coordinates{
(-5 , 2)
(0.382, 2)
(0.382, 2.5)
(5, 2.5)
		};
		\addplot[myblue, thick, mark=square*,mark size = 1.5pt] coordinates {
( -5.000000000000000 , 2.000000000000000)
( -4.814814814814815 , 2.000000000000000)
( -4.629629629629630 , 2.000000000000000)
( -4.444444444444445 , 2.000000000000000)
( -4.259259259259260 , 2.000000000000000)
( -4.074074074074074 , 2.000000000000000)
( -3.888888888888889 , 2.000000000000000)
( -3.703703703703704 , 2.000000000000000)
( -3.518518518518519 , 2.000000000000000)
( -3.333333333333333 , 2.000000000000000)
( -3.148148148148148 , 2.000000000000000)
( -2.962962962962963 , 2.000000000000000)
( -2.777777777777778 , 2.000000000000000)
( -2.592592592592593 , 2.000000000000000)
( -2.407407407407407 , 2.000000000000000)
( -2.222222222222222 , 2.000000000000000)
( -2.037037037037037 , 2.000000000000000)
( -1.851851851851852 , 2.000000000000000)
( -1.666666666666667 , 2.000000000000000)
( -1.481481481481482 , 2.000000000000000)
( -1.296296296296297 , 2.000000000000000)
( -1.111111111111111 , 2.000000000000000)
( -0.925925925925926 , 2.000000000000000)
( -0.740740740740741 , 2.000000000000000)
( -0.555555555555555 , 2.000000000000000)
( -0.370370370370370 , 2.000000000000000)
( -0.185185185185185 , 2.000000000000000)
( 0.000000000000000  , 2.045164917874771)
( 0.185185185185185  , 2.177545844417498)
( 0.370370370370370  , 2.339455807702551)
( 0.555555555555555  , 2.442895441250442)
( 0.740740740740740  , 2.483638588182234)
( 0.925925925925926  , 2.495227607013573)
( 1.111111111111111  , 2.497798806986666)
( 1.296296296296296  , 2.498197772586761)
( 1.481481481481481  , 2.498200895282965)
( 1.666666666666666  , 2.498206337272284)
( 1.851851851851851  , 2.498294890370750)
( 2.037037037037036  , 2.498470091688399)
( 2.222222222222221  , 2.498708005752043)
( 2.407407407407407  , 2.498973697536802)
( 2.592592592592593  , 2.499232817536277)
( 2.777777777777778  , 2.499459776158822)
( 2.962962962962963  , 2.499641171295620)
( 3.148148148148147  , 2.499774835645991)
( 3.333333333333332  , 2.499866308900697)
( 3.518518518518517  , 2.499924770044633)
( 3.703703703703702  , 2.499959816587693)
( 3.888888888888889  , 2.499979595303351)
( 4.074074074074074  , 2.499990135563950)
( 4.259259259259260  , 2.499995453536247)
( 4.444444444444445  , 2.499997999662693)
( 4.629629629629630  , 2.499999158802986)
( 4.814814814814815  , 2.499999661496402)
( 5.000000000000000  , 2.499999869510234)
			};
		\end{axis}
	\end{tikzpicture}
	}
	\hspace{2em}
	\subfloat[$v_G$ (Roe's method)]{
	\begin{tikzpicture}
		\begin{axis}
		\addplot[gray, thick, mark=none] coordinates{
(-5 , 1.5)
(0.382, 1.5)
(0.382, 1.2764)
(5, 1.2764)
		};
		\addplot[myblue, thick, mark=square*,mark size = 1.5pt] coordinates{
( -5.000000000000000 , 1.500000000000000)
( -4.814814814814815 , 1.500000000000000)
( -4.629629629629630 , 1.500000000000000)
( -4.444444444444445 , 1.500000000000000)
( -4.259259259259260 , 1.500000000000000)
( -4.074074074074074 , 1.500000000000000)
( -3.888888888888889 , 1.500000000000000)
( -3.703703703703704 , 1.500000000000000)
( -3.518518518518519 , 1.500000000000000)
( -3.333333333333333 , 1.500000000000000)
( -3.148148148148148 , 1.500000000000000)
( -2.962962962962963 , 1.500000000000000)
( -2.777777777777778 , 1.500000000000000)
( -2.592592592592593 , 1.500000000000000)
( -2.407407407407407 , 1.500000000000000)
( -2.222222222222222 , 1.500000000000000)
( -2.037037037037037 , 1.500000000000000)
( -1.851851851851852 , 1.500000000000000)
( -1.666666666666667 , 1.500000000000000)
( -1.481481481481482 , 1.500000000000000)
( -1.296296296296297 , 1.500000000000000)
( -1.111111111111111 , 1.500000000000000)
( -0.925925925925926 , 1.500000000000000)
( -0.740740740740741 , 1.500000000000000)
( -0.555555555555555 , 1.500000000000000)
( -0.370370370370370 , 1.500000000000000)
( -0.185185185185185 , 1.500000000000000)
( 0.000000000000000  , 1.477591371934628)
( 0.185185185185185  , 1.414395927671391)
( 0.370370370370370  , 1.342152063509673)
( 0.555555555555555  , 1.298659970764561)
( 0.740740740740740  , 1.281991695938680)
( 0.925925925925926  , 1.277207539561618)
( 1.111111111111111  , 1.276044582275619)
( 1.296296296296296  , 1.275768364446131)
( 1.481481481481481  , 1.275693736907137)
( 1.666666666666666  , 1.275684180296416)
( 1.851851851851851  , 1.275717992875401)
( 2.037037037037036  , 1.275787908400868)
( 2.222222222222221  , 1.275883097586719)
( 2.407407407407407  , 1.275989413824878)
( 2.592592592592593  , 1.276093092029172)
( 2.777777777777778  , 1.276183893978450)
( 2.962962962962963  , 1.276256461662854)
( 3.148148148148147  , 1.276309931734724)
( 3.333333333333332  , 1.276346522735855)
( 3.518518518518517  , 1.276369907778404)
( 3.703703703703702  , 1.276383926573167)
( 3.888888888888889  , 1.276391838107087)
( 4.074074074074074  , 1.276396054222662)
( 4.259259259259260  , 1.276398181413970)
( 4.444444444444445  , 1.276399199864992)
( 4.629629629629630  , 1.276399663521183)
( 4.814814814814815  , 1.276399864598558)
( 5.000000000000000  , 1.276399947804093)
			};
		\end{axis}
	\end{tikzpicture}
	}\\
	\subfloat[$m_L$ (Roe's method (squares), NT (triangles))]{
	\begin{tikzpicture}
		\begin{axis}[cycle list name=exotic]
		\addplot[gray, thick, mark=none] coordinates{
(-5,3)
(-2.2667, 3)
(-2.2667, 3.25)
(0.382, 3.25)
(0.382, 3.4995)
(3.5761, 3.4995)
(3.5761, 3)
(5, 3)
		};
		\addplot[myblue, thick, mark=square*,mark size = 1.5pt] coordinates{
( -4.907407407407407 , 3.000000009885070 )
( -4.722222222222222 , 3.000000043730698 )
( -4.537037037037037 , 3.000000183282335 )
( -4.351851851851852 , 3.000000730612745 )
( -4.166666666666667 , 3.000002772829557 )
( -3.981481481481482 , 3.000010005781617 )
( -3.796296296296296 , 3.000034275431503 )
( -3.611111111111111 , 3.000111272043823 )
( -3.425925925925926 , 3.000341729044733 )
( -3.240740740740741 , 3.000990867383230 )
( -3.055555555555556 , 3.002706207422170 )
( -2.870370370370371 , 3.006938761564575 )
( -2.685185185185185 , 3.016611566816831 )
( -2.500000000000000 , 3.036771086847558 )
( -2.314814814814815 , 3.073993786481731 )
( -2.129629629629630 , 3.131908106479002 )
( -1.944444444444445 , 3.202022776083672 )
( -1.759259259259260 , 3.259631971420604 )
( -1.574074074074074 , 3.280125597127556 )
( -1.388888888888889 , 3.264616931699845 )
( -1.203703703703704 , 3.240534150478267 )
( -1.018518518518519 , 3.233625904038871 )
( -0.833333333333334 , 3.244351925115326 )
( -0.648148148148149 , 3.254228631023792 )
( -0.462962962962963 , 3.252022811301083 )
( -0.277777777777778 , 3.245739732838930 )
( -0.092592592592593 , 3.262365338275981 )
( 0.092592592592593  , 3.324216420960761 )
( 0.277777777777778  , 3.404594580170750 )
( 0.462962962962963  , 3.462001048690390 )
( 0.648148148148148  , 3.487577214145956 )
( 0.833333333333333  , 3.496225968177875 )
( 1.018518518518518  , 3.496438482859451 )
( 1.203703703703703  , 3.496956006629335 )
( 1.388888888888888  , 3.497459070207306 )
( 1.574074074074074  , 3.496889750750211 )
( 1.759259259259259  , 3.496319489811058 )
( 1.944444444444444  , 3.495758022794626 )
( 2.129629629629629  , 3.492684521727786 )
( 2.314814814814814  , 3.488116540399742 )
( 2.499999999999999  , 3.478326801943377 )
( 2.685185185185185  , 3.462189009572917 )
( 2.870370370370370  , 3.436119569075935 )
( 3.055555555555555  , 3.398148002074802 )
( 3.240740740740740  , 3.347376857944446 )
( 3.425925925925925  , 3.286307216773318 )
( 3.611111111111110  , 3.220771044286967 )
( 3.796296296296295  , 3.158507946274389 )
( 3.981481481481482  , 3.106049864365748 )
( 4.166666666666667  , 3.066479113502207 )
( 4.351851851851852  , 3.039341147796306 )
( 4.537037037037037  , 3.022139470632296 )
( 4.722222222222222  , 3.011916543989265 )
( 4.907407407407407  , 3.006159003201680 )
			};
\addplot[myorange, thick, mark=triangle*] coordinates{
( -4.907407407407407 , 3.000072974038833 )
( -4.722222222222222 , 3.000208258321827 )
( -4.537037037037037 , 3.000418162919122 )
( -4.351851851851852 , 3.000761356138567 )
( -4.166666666666667 , 3.001348658497848 )
( -3.981481481481482 , 3.002353106899368 )
( -3.796296296296296 , 3.004033426668515 )
( -3.611111111111111 , 3.006771018795166 )
( -3.425925925925926 , 3.011114945333815 )
( -3.240740740740741 , 3.017825297312072 )
( -3.055555555555556 , 3.027899280852606 )
( -2.870370370370371 , 3.042554550397186 )
( -2.685185185185185 , 3.063132304597028 )
( -2.500000000000000 , 3.090889762107381 )
( -2.314814814814815 , 3.126458103422107 )
( -2.129629629629630 , 3.160525325603390 )
( -1.944444444444445 , 3.186224992835344 )
( -1.759259259259260 , 3.205147125184533 )
( -1.574074074074074 , 3.218986219518713 )
( -1.388888888888889 , 3.228982133874001 )
( -1.203703703703704 , 3.235848638598773 )
( -1.018518518518519 , 3.240186778193798 )
( -0.833333333333334 , 3.242899418714785 )
( -0.648148148148149 , 3.245058412746645 )
( -0.462962962962963 , 3.247387732500815 )
( -0.277777777777778 , 3.250960199613859 )
( -0.092592592592593 , 3.259050165112657 )
( 0.092592592592593  , 3.281798010566062 )
( 0.277777777777778  , 3.332984661900177 )
( 0.462962962962963  , 3.406964314538778 )
( 0.648148148148148  , 3.457052257037620 )
( 0.833333333333333  , 3.482178731883007 )
( 1.018518518518518  , 3.491642701221142 )
( 1.203703703703703  , 3.493915088619490 )
( 1.388888888888888  , 3.494162657841982 )
( 1.574074074074074  , 3.493959423897690 )
( 1.759259259259259  , 3.492979122477380 )
( 1.944444444444444  , 3.490807284390911 )
( 2.129629629629629  , 3.487352060560788 )
( 2.314814814814814  , 3.482545733351796 )
( 2.499999999999999  , 3.475811431655464 )
( 2.685185185185185  , 3.465578123275185 )
( 2.870370370370370  , 3.449099835341245 )
( 3.055555555555555  , 3.422607750860758 )
( 3.240740740740740  , 3.381565734521933 )
( 3.425925925925925  , 3.320920575930807 )
( 3.611111111111110  , 3.238736877860710 )
( 3.796296296296295  , 3.161537729640331 )
( 3.981481481481482  , 3.105513386324442 )
( 4.166666666666667  , 3.066566618409278 )
( 4.351851851851852  , 3.040281580438868 )
( 4.537037037037037  , 3.023144524676038 )
( 4.722222222222222  , 3.013154773211277 )
( 4.907407407407407  , 3.009181668494195 )
            };
		\end{axis}
	\end{tikzpicture}
	}
	\hspace{2em}
	\subfloat[$v_L$ (Roe's method (squares), NT (triangles))]{
	\begin{tikzpicture}
		\begin{axis}
		\addplot[gray, thick, mark=none] coordinates{
(-5,1)
(-2.2667, 1)
(-2.2667, 0.7487)
(0.382, 0.7487)
(0.382, 0.7226)
(3.5761, 0.7226)
(3.5761, 0.2475)
(5, 0.2475)
		};

		\addplot[myblue, thick,, mark=square*,mark size = 1.5pt] coordinates {
( -4.907407407407407 , 0.999999993071306 )
( -4.722222222222222 , 0.999999971381150 )
( -4.537037037037037 , 0.999999878889578 )
( -4.351851851851852 , 0.999999509670708 )
( -4.166666666666667 , 0.999998107979462 )
( -3.981481481481482 , 0.999993054548125 )
( -3.796296296296296 , 0.999975783899066 )
( -3.611111111111111 , 0.999919942911558 )
( -3.425925925925926 , 0.999749503337827 )
( -3.240740740740741 , 0.999259633076870 )
( -3.055555555555556 , 0.997937903992877 )
( -2.870370370370371 , 0.994604986761691 )
( -2.685185185185185 , 0.986808466240963 )
( -2.500000000000000 , 0.970116422750156 )
( -2.314814814814815 , 0.938193294711145 )
( -2.129629629629630 , 0.885819241539553 )
( -1.944444444444445 , 0.816599461298749 )
( -1.759259259259260 , 0.750382561433290 )
( -1.574074074074074 , 0.716117932334631 )
( -1.388888888888889 , 0.724958365321765 )
( -1.203703703703704 , 0.753320544109907 )
( -1.018518518518519 , 0.767485055128086 )
( -0.833333333333334 , 0.758726007302053 )
( -0.648148148148149 , 0.744866953492534 )
( -0.462962962962963 , 0.743941591004848 )
( -0.277777777777778 , 0.752087845039504 )
( -0.092592592592593 , 0.770206733327287 )
( 0.092592592592593  , 0.777836410128016 )
( 0.277777777777778  , 0.765730079563770 )
( 0.462962962962963  , 0.743096601088781 )
( 0.648148148148148  , 0.730569284287633 )
( 0.833333333333333  , 0.723040241660202 )
( 1.018518518518518  , 0.723936869213881 )
( 1.203703703703703  , 0.723557840193135 )
( 1.388888888888888  , 0.722772650022861 )
( 1.574074074074074  , 0.722367761432597 )
( 1.759259259259259  , 0.723070264317537 )
( 1.944444444444444  , 0.720249650427510 )
( 2.129629629629629  , 0.719278771314757 )
( 2.314814814814814  , 0.713172343044626 )
( 2.499999999999999  , 0.704637236472224 )
( 2.685185185185185  , 0.688464961390606 )
( 2.870370370370370  , 0.663623325548520 )
( 3.055555555555555  , 0.626957316746916 )
( 3.240740740740740  , 0.578238496908338 )
( 3.425925925925925  , 0.519610977791021 )
( 3.611111111111110  , 0.456858018512725 )
( 3.796296296296295  , 0.397401278438848 )
( 3.981481481481482  , 0.347509017290880 )
( 4.166666666666667  , 0.310032582727976 )
( 4.351851851851852  , 0.284431745791069 )
( 4.537037037037037  , 0.268254172952365 )
( 4.722222222222222  , 0.258660788326197 )
( 4.907407407407407  , 0.253265760625346 )
			};
\addplot[myorange, thick,, mark=triangle*] coordinates {
( -4.907407407407407 , 0.999927021367714 )
( -4.722222222222222 , 0.999791726497481 )
( -4.537037037037037 , 0.999581793216027 )
( -4.351851851851852 , 0.999238512145331 )
( -4.166666666666667 , 0.998650949899532 )
( -3.981481481481482 , 0.997645779244333 )
( -3.796296296296296 , 0.995963563556076 )
( -3.611111111111111 , 0.993221249985417 )
( -3.425925925925926 , 0.988866182119251 )
( -3.240740740740741 , 0.982131059415653 )
( -3.055555555555556 , 0.972005764045186 )
( -2.870370370370371 , 0.957253961699721 )
( -2.685185185185185 , 0.936523586432272 )
( -2.500000000000000 , 0.908644766416932 )
( -2.314814814814815 , 0.873407497535985 )
( -2.129629629629630 , 0.838779245063499 )
( -1.944444444444445 , 0.812197999264030 )
( -1.759259259259260 , 0.792848677133643 )
( -1.574074074074074 , 0.779366433572130 )
( -1.388888888888889 , 0.770190951973509 )
( -1.203703703703704 , 0.763754265274246 )
( -1.018518518518519 , 0.758915296918038 )
( -0.833333333333334 , 0.755298865670371 )
( -0.648148148148149 , 0.753067253851030 )
( -0.462962962962963 , 0.752143754944211 )
( -0.277777777777778 , 0.752413044570557 )
( -0.092592592592593 , 0.757935930196915 )
( 0.092592592592593  , 0.772702539616753 )
( 0.277777777777778  , 0.774701997227183 )
( 0.462962962962963  , 0.760309119286944 )
( 0.648148148148148  , 0.747248227445285 )
( 0.833333333333333  , 0.734559547296878 )
( 1.018518518518518  , 0.725576665043203 )
( 1.203703703703703  , 0.721369323927650 )
( 1.388888888888888  , 0.720083199103847 )
( 1.574074074074074  , 0.718919768315785 )
( 1.759259259259259  , 0.717559272279783 )
( 1.944444444444444  , 0.716127087620434 )
( 2.129629629629629  , 0.713980769864296 )
( 2.314814814814814  , 0.710005199645089 )
( 2.499999999999999  , 0.702973276494512 )
( 2.685185185185185  , 0.691631406681531 )
( 2.870370370370370  , 0.674296650920172 )
( 3.055555555555555  , 0.648178296320428 )
( 3.240740740740740  , 0.608990840891827 )
( 3.425925925925925  , 0.551212768451571 )
( 3.611111111111110  , 0.472913447952317 )
( 3.796296296296295  , 0.400881295554801 )
( 3.981481481481482  , 0.347534262924103 )
( 4.166666666666667  , 0.310418106416803 )
( 4.351851851851852  , 0.285472739741759 )
( 4.537037037037037  , 0.269272928117229 )
( 4.722222222222222  , 0.259852746020838 )
( 4.907407407407407  , 0.256112031498598 )
            };
		\end{axis}
	\end{tikzpicture}
	}
    \caption{Exact and numerical solutions for the all-shock Riemann problem.}
    \label{fig: Riemann problem 1}
\end{figure}
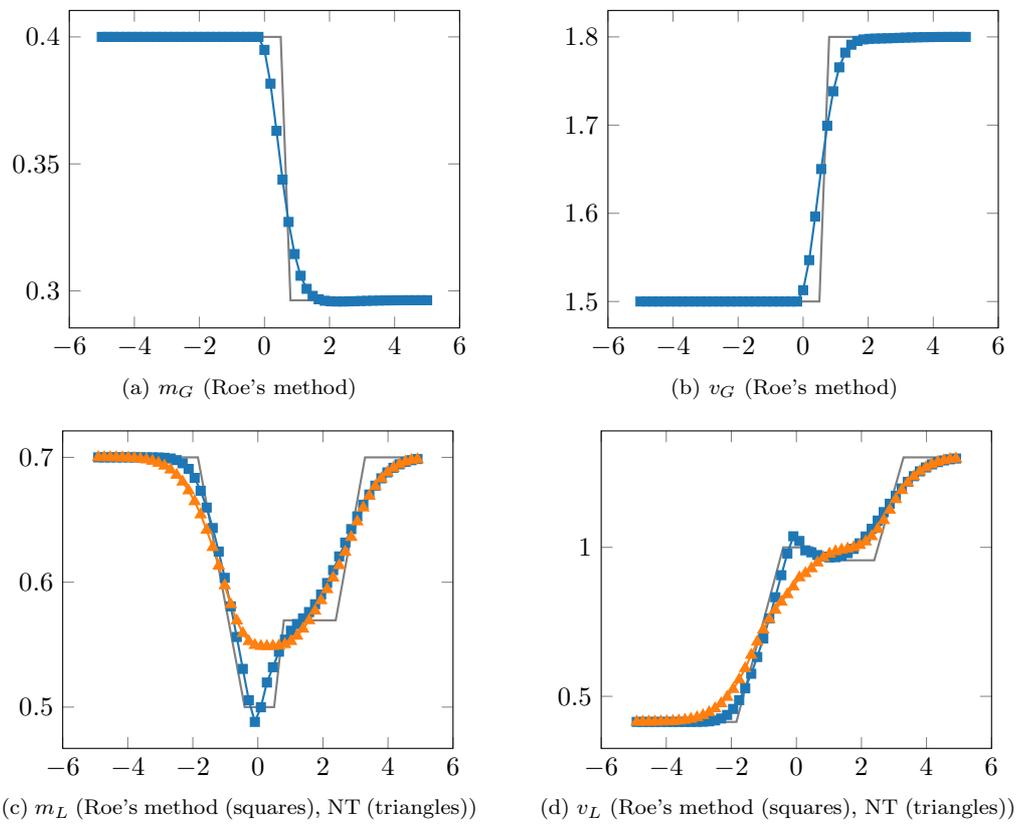
\begin{figure}
    \centering
    \subfloat[$m_G$ (Roe's method)]{
    \begin{tikzpicture}
        \begin{axis}
        \addplot[gray, thick, mark=none] coordinates{
(-5 , 0.4)
(0.5, 0.4)
(0.8, 0.2963)
(5, 0.2963)
        };
        \addplot[myblue, thick, mark=square*,mark size = 1.5pt] coordinates {
( -5.000000000000000 , 0.400000000000000 )
( -4.814814814814815 , 0.400000000000000 )
( -4.629629629629630 , 0.400000000000000 )
( -4.444444444444445 , 0.400000000000000 )
( -4.259259259259260 , 0.400000000000000 )
( -4.074074074074074 , 0.400000000000000 )
( -3.888888888888889 , 0.400000000000000 )
( -3.703703703703704 , 0.400000000000000 )
( -3.518518518518519 , 0.400000000000000 )
( -3.333333333333333 , 0.400000000000000 )
( -3.148148148148148 , 0.400000000000000 )
( -2.962962962962963 , 0.400000000000000 )
( -2.777777777777778 , 0.400000000000000 )
( -2.592592592592593 , 0.400000000000000 )
( -2.407407407407407 , 0.400000000000000 )
( -2.222222222222222 , 0.400000000000000 )
( -2.037037037037037 , 0.400000000000000 )
( -1.851851851851852 , 0.400000000000000 )
( -1.666666666666667 , 0.400000000000000 )
( -1.481481481481482 , 0.400000000000000 )
( -1.296296296296297 , 0.400000000000000 )
( -1.111111111111111 , 0.400000000000000 )
( -0.925925925925926 , 0.400000000000000 )
( -0.740740740740741 , 0.400000000000000 )
( -0.555555555555555 , 0.400000000000000 )
( -0.370370370370370 , 0.400000000000000 )
( -0.185185185185185 , 0.400000000000000 )
( 0.000000000000000  , 0.394891788455086 )
( 0.185185185185185  , 0.381594077019649 )
( 0.370370370370370  , 0.363027569872844 )
( 0.555555555555555  , 0.343810534687973 )
( 0.740740740740740  , 0.327159203103432 )
( 0.925925925925926  , 0.314509214185762 )
( 1.111111111111111  , 0.305962925168721 )
( 1.296296296296296  , 0.300814265992474 )
( 1.481481481481481  , 0.298043877806814 )
( 1.666666666666666  , 0.296705785067775 )
( 1.851851851851851  , 0.296121876520714 )
( 2.037037037037036  , 0.295894597575501 )
( 2.222222222222221  , 0.295826764868815 )
( 2.407407407407407  , 0.295832378579322 )
( 2.592592592592593  , 0.295875509868595 )
( 2.777777777777778  , 0.295938774856007 )
( 2.962962962962963  , 0.296010749705135 )
( 3.148148148148147  , 0.296082520388561 )
( 3.333333333333332  , 0.296147494877133 )
( 3.518518518518517  , 0.296201732957152 )
( 3.703703703703702  , 0.296243857236060 )
( 3.888888888888889  , 0.296274482850464 )
( 4.074074074074074  , 0.296295421504435 )
( 4.259259259259260  , 0.296308932839472 )
( 4.444444444444445  , 0.296317186156581 )
( 4.629629629629630  , 0.296321970641154 )
( 4.814814814814815  , 0.296324608576157 )
( 5.000000000000000  , 0.296325994517794 )
            };
        \end{axis}
    \end{tikzpicture}
    }
    \hspace{2em}
    \subfloat[$v_G$ (Roe's method)]{
    \begin{tikzpicture}
        \begin{axis}
        \addplot[gray, thick, mark=none] coordinates{
(-5 , 1.5)
(0.5, 1.5)
(0.8, 1.8)
(5,  1.8)
        };
        \addplot[myblue, thick, mark=square*, mark size=1.5pt] coordinates{
( -5.000000000000000 , 1.500000000000000)
( -4.814814814814815 , 1.500000000000000)
( -4.629629629629630 , 1.500000000000000)
( -4.444444444444445 , 1.500000000000000)
( -4.259259259259260 , 1.500000000000000)
( -4.074074074074074 , 1.500000000000000)
( -3.888888888888889 , 1.500000000000000)
( -3.703703703703704 , 1.500000000000000)
( -3.518518518518519 , 1.500000000000000)
( -3.333333333333333 , 1.500000000000000)
( -3.148148148148148 , 1.500000000000000)
( -2.962962962962963 , 1.500000000000000)
( -2.777777777777778 , 1.500000000000000)
( -2.592592592592593 , 1.500000000000000)
( -2.407407407407407 , 1.500000000000000)
( -2.222222222222222 , 1.500000000000000)
( -2.037037037037037 , 1.500000000000000)
( -1.851851851851852 , 1.500000000000000)
( -1.666666666666667 , 1.500000000000000)
( -1.481481481481482 , 1.500000000000000)
( -1.296296296296297 , 1.500000000000000)
( -1.111111111111111 , 1.500000000000000)
( -0.925925925925926 , 1.500000000000000)
( -0.740740740740741 , 1.500000000000000)
( -0.555555555555555 , 1.500000000000000)
( -0.370370370370370 , 1.500000000000000)
( -0.185185185185185 , 1.500000000000000)
( 0.000000000000000  , 1.512825117141126)
( 0.185185185185185  , 1.546895530677711)
( 0.370370370370370  , 1.596362750796530)
( 0.555555555555555  , 1.650205715442588)
( 0.740740740740740  , 1.699317056771083)
( 0.925925925925926  , 1.738312903340276)
( 1.111111111111111  , 1.765513596242052)
( 1.296296296296296  , 1.782179851298145)
( 1.481481481481481  , 1.791143395383232)
( 1.666666666666666  , 1.795378754043306)
( 1.851851851851851  , 1.797144747950327)
( 2.037037037037036  , 1.797814005850808)
( 2.222222222222221  , 1.798086060465798)
( 2.407407407407407  , 1.798264842592880)
( 2.592592592592593  , 1.798457469242170)
( 2.777777777777778  , 1.798683945366389)
( 2.962962962962963  , 1.798930285144905)
( 3.148148148148147  , 1.799173425434375)
( 3.333333333333332  , 1.799392977775474)
( 3.518518518518517  , 1.799576114361625)
( 3.703703703703702  , 1.799718309233163)
( 3.888888888888889  , 1.799821675798130)
( 4.074074074074074  , 1.799892342271974)
( 4.259259259259260  , 1.799937940251231)
( 4.444444444444445  , 1.799965792885708)
( 4.629629629629630  , 1.799981938994107)
( 4.814814814814815  , 1.799990841131789)
( 5.000000000000000  , 1.799995518203375)
            };
        \end{axis}
    \end{tikzpicture}
    }\\
    \subfloat[$m_L$ (Roe's method (squares), NT (triangles))]{
    \begin{tikzpicture}
        \begin{axis}
        \addplot[gray, thick, mark=none] coordinates{
(-5,0.7)
(-1.8441, 0.7)
(-0.4053, 0.5)
(0.5, 0.5)
(0.8, 0.5695)
(2.3936, 0.5695)
(3.2941, 0.7)
(5, 0.7)
        };
        \addplot[myblue, thick, mark=square*,mark size = 1.5pt] coordinates{
( -4.907407407407407 , 0.699999999897991 )
( -4.722222222222222 , 0.699999999459750 )
( -4.537037037037037 , 0.699999997307256 )
( -4.351851851851852 , 0.699999987272711 )
( -4.166666666666667 , 0.699999942985680 )
( -3.981481481481482 , 0.699999758488035 )
( -3.796296296296296 , 0.699999035199875 )
( -3.611111111111111 , 0.699996376274138 )
( -3.425925925925926 , 0.699987248125389 )
( -3.240740740740741 , 0.699958126106066 )
( -3.055555555555556 , 0.699872282788318 )
( -2.870370370370371 , 0.699640073841330 )
( -2.685185185185185 , 0.699068096034819 )
( -2.500000000000000 , 0.697795527050968 )
( -2.314814814814815 , 0.695256715198301 )
( -2.129629629629630 , 0.690735060199624 )
( -1.944444444444445 , 0.683541247930992 )
( -1.759259259259260 , 0.673241852989937 )
( -1.574074074074074 , 0.659783261419795 )
( -1.388888888888889 , 0.643424220281659 )
( -1.203703703703704 , 0.624551480199288 )
( -1.018518518518519 , 0.603524730230752 )
( -0.833333333333334 , 0.580630654140399 )
( -0.648148148148149 , 0.556148578316056 )
( -0.462962962962963 , 0.530580828924522 )
( -0.277777777777778 , 0.505490685746027 )
( -0.092592592592593 , 0.488123601370701 )
( 0.092592592592593  , 0.499867723688391 )
( 0.277777777777778  , 0.519702600201068 )
( 0.462962962962963  , 0.531946037409993 )
( 0.648148148148148  , 0.544389751206059 )
( 0.833333333333333  , 0.554052357960400 )
( 1.018518518518518  , 0.561168692969700 )
( 1.203703703703703  , 0.566462452007044 )
( 1.388888888888888  , 0.571075311945268 )
( 1.574074074074074  , 0.576120260053971 )
( 1.759259259259259  , 0.582374364376245 )
( 1.944444444444444  , 0.590158744479094 )
( 2.129629629629629  , 0.599388907999776 )
( 2.314814814814814  , 0.609713474168651 )
( 2.499999999999999  , 0.620661074759824 )
( 2.685185185185185  , 0.631750562108564 )
( 2.870370370370370  , 0.642554990175370 )
( 3.055555555555555  , 0.652729337526077 )
( 3.240740740740740  , 0.662016628400172 )
( 3.425925925925925  , 0.670244129252413 )
( 3.611111111111110  , 0.677316398414486 )
( 3.796296296296295  , 0.683207822522975 )
( 3.981481481481482  , 0.687954609436218 )
( 4.166666666666667  , 0.691645030581395 )
( 4.351851851851852  , 0.694406847710942 )
( 4.537037037037037  , 0.696391979195726 )
( 4.722222222222222  , 0.697759926405320 )
( 4.907407407407407  , 0.698662472886416 )
            };
\addplot[myorange, thick, mark=triangle*]coordinates{
( -4.907407407407407 , 0.699986918019005 )
( -4.722222222222222 , 0.699962356329051 )
( -4.537037037037037 , 0.699924129438964 )
( -4.351851851851852 , 0.699862155701842 )
( -4.166666666666667 , 0.699757513758528 )
( -3.981481481481482 , 0.699581201358662 )
( -3.796296296296296 , 0.699291211927927 )
( -3.611111111111111 , 0.698828031605735 )
( -3.425925925925926 , 0.698109530739093 )
( -3.240740740740741 , 0.697026510302545 )
( -3.055555555555556 , 0.695440125513583 )
( -2.870370370370371 , 0.693182220490863 )
( -2.685185185185185 , 0.690059255123995 )
( -2.500000000000000 , 0.685859930563231 )
( -2.314814814814815 , 0.680365975990619 )
( -2.129629629629630 , 0.673365498883455 )
( -1.944444444444445 , 0.664670552156147 )
( -1.759259259259260 , 0.654149292566377 )
( -1.574074074074074 , 0.641805039922240 )
( -1.388888888888889 , 0.627958642440335 )
( -1.203703703703704 , 0.613125959309107 )
( -1.018518518518519 , 0.597369270955206 )
( -0.833333333333334 , 0.582181397048455 )
( -0.648148148148149 , 0.569173853778013 )
( -0.462962962962963 , 0.559138592163965 )
( -0.277777777777778 , 0.552613304608110 )
( -0.092592592592593 , 0.549541979200944 )
( 0.092592592592593  , 0.548742414387413 )
( 0.277777777777778  , 0.548707293429424 )
( 0.462962962962963  , 0.548783196855270 )
( 0.648148148148148  , 0.549079230204193 )
( 0.833333333333333  , 0.550238114344943 )
( 1.018518518518518  , 0.552702173637920 )
( 1.203703703703703  , 0.556948478448789 )
( 1.388888888888888  , 0.562598789505654 )
( 1.574074074074074  , 0.569356383349610 )
( 1.759259259259259  , 0.577024162602326 )
( 1.944444444444444  , 0.585367755312788 )
( 2.129629629629629  , 0.594219045807790 )
( 2.314814814814814  , 0.603535731273433 )
( 2.499999999999999  , 0.613408020811326 )
( 2.685185185185185  , 0.624108392844739 )
( 2.870370370370370  , 0.636039557992451 )
( 3.055555555555555  , 0.648533697819064 )
( 3.240740740740740  , 0.659623974197569 )
( 3.425925925925925  , 0.669032165810326 )
( 3.611111111111110  , 0.676769325437774 )
( 3.796296296296295  , 0.682947713812419 )
( 3.981481481481482  , 0.687761082713753 )
( 4.166666666666667  , 0.691460828600322 )
( 4.351851851851852  , 0.694290414197253 )
( 4.537037037037037  , 0.696387731576563 )
( 4.722222222222222  , 0.697741863511872 )
( 4.907407407407407  , 0.698314855803901 )
            };
        \end{axis}
    \end{tikzpicture}
    }
    \hspace{2em}
    \subfloat[$v_L$ (Roe's method (squares), NT (triangles))]{
    \begin{tikzpicture}
        \begin{axis}
        \addplot[gray, thick, mark=none] coordinates{
(-5,0.4141)
(-1.8441, 0.4141)
(-0.4053, 1.)
(0.5, 1.)
(0.8, 0.9566)
(2.3936, 0.9566)
(3.2941, 1.3021)
(5, 1.3021)
        };

        \addplot[myblue, thick, mark=square*,mark size = 1.5pt] coordinates {
( -4.907407407407407 , 0.414084440173840 )
( -4.722222222222222 , 0.414084440849427 )
( -4.537037037037037 , 0.414084444322000 )
( -4.351851851851852 , 0.414084460951031 )
( -4.166666666666667 , 0.414084536394083 )
( -3.981481481481482 , 0.414084859911678 )
( -3.796296296296296 , 0.414086167496529 )
( -3.611111111111111 , 0.414091132246642 )
( -3.425925925925926 , 0.414108772459009 )
( -3.240740740740741 , 0.414167159109700 )
( -3.055555555555556 , 0.414346212030403 )
( -2.870370370370371 , 0.414851717731910 )
( -2.685185185185185 , 0.416155802152875 )
( -2.500000000000000 , 0.419204443617220 )
( -2.314814814814815 , 0.425608363004768 )
( -2.129629629629630 , 0.437608890202844 )
( -1.944444444444445 , 0.457597804100991 )
( -1.759259259259260 , 0.487270185658196 )
( -1.574074074074074 , 0.526954587449220 )
( -1.388888888888889 , 0.575678914785555 )
( -1.203703703703704 , 0.631876126116951 )
( -1.018518518518519 , 0.694093596588293 )
( -0.833333333333334 , 0.761262798780839 )
( -0.648148148148149 , 0.832523752132499 )
( -0.462962962962963 , 0.906552158685253 )
( -0.277777777777778 , 0.979192330127649 )
( -0.092592592592593 , 1.037129869292003 )
( 0.092592592592593  , 1.021518848580388 )
( 0.277777777777778  , 0.989929372460168 )
( 0.462962962962963  , 0.983245491551986 )
( 0.648148148148148  , 0.974220918612202 )
( 0.833333333333333  , 0.968759293408163 )
( 1.018518518518518  , 0.966233350575617 )
( 1.203703703703703  , 0.967002877707399 )
( 1.388888888888888  , 0.971639928284485 )
( 1.574074074074074  , 0.980787616675882 )
( 1.759259259259259  , 0.994817869989727 )
( 1.944444444444444  , 1.013597081645845 )
( 2.129629629629629  , 1.036461636975632 )
( 2.314814814814814  , 1.062369703452447 )
( 2.499999999999999  , 1.090114415328691 )
( 2.685185185185185  , 1.118502760141576 )
( 2.870370370370370  , 1.146465246348449 )
( 3.055555555555555  , 1.173105952840905 )
( 3.240740740740740  , 1.197717228237514 )
( 3.425925925925925  , 1.219779664527747 )
( 3.611111111111110  , 1.238958662287859 )
( 3.796296296296295  , 1.255100349901497 )
( 3.981481481481482  , 1.268223913132739 )
( 4.166666666666667  , 1.278505119719094 )
( 4.351851851851852  , 1.286247112569055 )
( 4.537037037037037  , 1.291838721825616 )
( 4.722222222222222  , 1.295705684325429 )
( 4.907407407407407  , 1.298263622750766 )
            };
\addplot[myorange, thick,, mark=triangle*] coordinates {
( -4.907407407407407 , 0.414126637489436 )
( -4.722222222222222 , 0.414205859269690 )
( -4.537037037037037 , 0.414329138766623 )
( -4.351851851851852 , 0.414528947084076 )
( -4.166666666666667 , 0.414866182165138 )
( -3.981481481481482 , 0.415434042375866 )
( -3.796296296296296 , 0.416367165823853 )
( -3.611111111111111 , 0.417855548410808 )
( -3.425925925925926 , 0.420159849101434 )
( -3.240740740740741 , 0.423623673460124 )
( -3.055555555555556 , 0.428678517665576 )
( -2.870370370370371 , 0.435837797443893 )
( -2.685185185185185 , 0.445677969986401 )
( -2.500000000000000 , 0.458807583801955 )
( -2.314814814814815 , 0.475829234237890 )
( -2.129629629629630 , 0.497304850135664 )
( -1.944444444444445 , 0.523743449989550 )
( -1.759259259259260 , 0.555651792961985 )
( -1.574074074074074 , 0.593745500543924 )
( -1.388888888888889 , 0.638446900524402 )
( -1.203703703703704 , 0.682909739266685 )
( -1.018518518518519 , 0.723319864205589 )
( -0.833333333333334 , 0.758634992346958 )
( -0.648148148148149 , 0.789466752516584 )
( -0.462962962962963 , 0.816894056150158 )
( -0.277777777777778 , 0.841912315016880 )
( -0.092592592592593 , 0.868123106094303 )
( 0.092592592592593  , 0.898475770459551 )
( 0.277777777777778  , 0.912581913645766 )
( 0.462962962962963  , 0.928843111134715 )
( 0.648148148148148  , 0.948368965529035 )
( 0.833333333333333  , 0.965550304148088 )
( 1.018518518518518  , 0.978079580650980 )
( 1.203703703703703  , 0.986163815948729 )
( 1.388888888888888  , 0.991409886099163 )
( 1.574074074074074  , 0.995326730457977 )
( 1.759259259259259  , 0.999837073125190 )
( 1.944444444444444  , 1.007120649014315 )
( 2.129629629629629  , 1.019036567026236 )
( 2.314814814814814  , 1.036795886035927 )
( 2.499999999999999  , 1.060823812300211 )
( 2.685185185185185  , 1.090720220132034 )
( 2.870370370370370  , 1.125503824805040 )
( 3.055555555555555  , 1.160985325690216 )
( 3.240740740740740  , 1.191325349871319 )
( 3.425925925925925  , 1.216757528441246 )
( 3.611111111111110  , 1.237710222371474 )
( 3.796296296296295  , 1.254559614481796 )
( 3.981481481481482  , 1.267799127771754 )
( 4.166666666666667  , 1.278061685327728 )
( 4.351851851851852  , 1.285966993841690 )
( 4.537037037037037  , 1.291856139917340 )
( 4.722222222222222  , 1.295669882912694 )
( 4.907407407407407  , 1.297286700211852 )
            };
        \end{axis}
    \end{tikzpicture}
    }
    \caption{Exact and numerical solutions for the all-rarefaction Riemann problem.}
    \label{fig: Riemann problem 2}
\end{figure}

% Roe's method, being an approximate Riemann solver, resolves the shocks better than the Nessyahu--Tadmor scheme which here suffers from numerical viscosity. The rarefaction wave that is present in the exact solution, on the other hand, is resolved slightly worse by Roe's method.

Tables \ref{tbl: L1 errors all-shock} and \ref{tbl: L1 errors all-rarefaction} show the relative $L^1$ errors between each of the two numerical approximations and the exact solution of the liquid mass and liquid velocity in Experiment 1 and 2 respectively for several values of the spatial grid size. We observe that the Roe scheme as an approximate Riemann solver has an overall smaller error in all instances.
\begin{table}
% \centering
% \renewcommand{\arraystretch}{1.3}
\centering
\subfloat[Liquid mass.]{
\begin{tabular}{rcc}
  \toprule
  \multicolumn{1}{c}{$N$} & rel. $L^1$ error (Roe) & rel. $L^1$ error (NT)\\
  \midrule
     $ 16$&  $2.93$ & $3.52$ \\
     $ 32$&  $1.81$ & $2.24$ \\
     $ 64$&  $1.09$ & $1.31$ \\
     $128$&  $0.65$ & $0.72$ \\
     $256$&  $0.37$ & $0.39$ \\
  \bottomrule
\end{tabular}
}
\subfloat[Liquid velocity.]{
\begin{tabular}{rcc}
  \toprule
  \multicolumn{1}{c}{$N$} & rel. $L^1$ error (Roe) & rel. $L^1$ error (NT)\\
  \midrule
     $ 16$&  $10.55$ & $12.26$ \\
     $ 32$&  $\phantom{1}6.46$ & $\phantom{1}8.72$ \\
     $ 64$&  $\phantom{1}3.90$ & $\phantom{1}5.13$ \\
     $128$&  $\phantom{1}2.34$ & $\phantom{1}2.81$ \\
     $256$&  $\phantom{1}1.34$ & $\phantom{1}1.48$ \\
  \bottomrule
\end{tabular}  
}
\caption{$L^1$ error between numerical approximations and exact solution of Experiment 1. ($\lambda = 4$)}
\label{tbl: L1 errors all-shock}
\end{table}
\begin{table}
% \centering
% \renewcommand{\arraystretch}{1.3}
\centering
\subfloat[Liquid mass.]{
\begin{tabular}{rcc}
  \toprule
  \multicolumn{1}{c}{$N$} & rel. $L^1$ error (Roe) & rel. $L^1$ error (NT)\\
  \midrule
     $ 16$&  $3.67$ & $7.53$ \\
     $ 32$&  $2.29$ & $4.68$ \\
     $ 64$&  $1.57$ & $2.77$ \\
     $128$&  $1.10$ & $1.59$ \\
     $256$&  $0.80$ & $0.89$ \\
  \bottomrule
\end{tabular}
}
\subfloat[Liquid velocity.]{
\begin{tabular}{rcc}
  \toprule
  \multicolumn{1}{c}{$N$} & rel. $L^1$ error (Roe) & rel. $L^1$ error (NT)\\
  \midrule
     $ 16$&  $6.78$ & $11.16$ \\
     $ 32$&  $4.42$ & $\phantom{1}8.44$ \\
     $ 64$&  $2.92$ & $\phantom{1}5.76$ \\
     $128$&  $1.99$ & $\phantom{1}3.09$ \\
     $256$&  $1.42$ & $\phantom{1}1.70$ \\
  \bottomrule
\end{tabular}  
}
\caption{$L^1$ error between numerical approximations and exact solution of Experiment 2. ($\lambda = 4$)}
\label{tbl: L1 errors all-rarefaction}
\end{table}

\section{Conclusion}
We have considered a model for two-phase flow in pipes where we assumed hydrostatic balance, that the gas is ideal, and that the liquid incompressible.
These assumptions allow us to partially decouple the $4\times 4$ system. Therefore, we can numerically treat the mass and momentum equation of the gas phase independently of the liquid phase. Staggering the numerical grid for the two subsystems with respect to each other further simplifies the numerical treatment since the approximation of the gas phase is constant at the cell interfaces of the approximation of the liquid phase.

As a starting point, we proposed to use the Roe scheme to numerically solve the gas subsystem and either the Roe scheme or the nonstaggered, second-order Nessyahu--Tadmor scheme for the liquid subsystem. In order to compare our numerical methods to exact solutions we further provided two classes of Riemann problems admitting only shocks respectively rarefactions. These classes can readily be modified to contain more than one shock, more than one rarefaction, or both shocks and rarefactions.

In our numerical experiments we compared numerical solutions of the liquid phase generated by the Roe scheme to the second-order Nessyahu--Tadmor scheme for certain Riemann problems.

% \subsection{Other initial data}

\section*{Acknowledgements}
We would like to thank Trygve Karper and Gunnar Staff from Schlumberger for helpful discussions on the topic.

% \newpage

\bibliographystyle{siam}
% \bibliography{literature}

\end{document}